\documentclass{amsart}

\input xy
\xyoption{all}

\usepackage[all]{xy}
\usepackage{tikz}
\usepackage{tikz-cd}
\usepackage{quiver}

\usepackage{enumitem}
\usepackage{geometry}
\usepackage{amsmath}
\usepackage{amssymb}
\usepackage{hyperref}
\usepackage{cite}
\usepackage[final]{showkeys}
\usepackage{hyperref}
\usepackage{setspace}
\usepackage{amsthm}  
\usepackage[pagewise,displaymath,mathlines]{lineno}

\newtheorem{same}{This should never appear}[section]
\newtheorem{defin}[same]{Definition}

\newtheorem{remark}[same]{Remark}

\newtheorem{example}[same]{Example}
\newtheorem{lemma}[same]{Lemma}
\newtheorem{fact}[same]{Fact}
\newtheorem{question}[same]{Question}

\newtheorem{prop}[same]{Proposition}

\newtheorem{conj}[same]{Conjecture}

\newbox\noforkbox \newdimen\forklinewidth
\setbox0\hbox{$\textstyle\smile$}
\setbox1\hbox to \wd0{\hfil\vrule width \forklinewidth depth-2pt
 height 10pt \hfil}
\setbox\noforkbox\hbox{\lower 2pt\box1\lower 2pt\box0\relax}
\def\unionstick{\mathop{\copy\noforkbox}\limits}

\def\nonfork_#1{\unionstick_{\textstyle #1}}

\setbox0\hbox{$\textstyle\smile$}
\setbox1\hbox to \wd0{\hfil{\sl /\/}\hfil}
\setbox2\hbox to \wd0{\hfil\vrule height 10pt depth -2pt width
              \forklinewidth\hfil}
\newbox\doesforkbox
\setbox\doesforkbox\hbox{\lower 2pt\box1 \lower 2pt\box2\lower2pt\box0\relax}
\def\nunionstick{\mathop{\copy\doesforkbox}\limits}

\def\fork_#1{\nunionstick_{\textstyle #1}}

\newcommand{\LS}{\text{LS}}

\newcommand{\ba}{\bold{a}}
\newcommand{\bb}{\bold{b}}
\newcommand{\bc}{\bold{c}}

\newcommand{\bm}{\bold{m}}

\newcommand{\cK}{\mathcal{K}}
\newcommand{\bK}{\mathbb{K}}
\newcommand{\bx}{\bold{x}}
\newcommand{\by}{\bold{y}}
\newcommand{\bz}{\bold{z}}

\newcommand{\cf}{\text{cf }}
\newcommand{\rest}{\upharpoonright}
\newcommand{\id}{\textrm{id}}

\newcommand{\gtp}{\text{gtp}}
\newcommand{\ZFC}{\text{ZFC}}
\newcommand{\cl}{\text{cl}}
\newcommand{\qftp}{\text{qftp}}

\renewcommand{\S}{\text{S}}

 \newcommand{\footnotei}[1]{}

\newcommand{\seq}[1]{\langle #1 \rangle}

\newcommand{\bL}{\mathbb{L}}
\newcommand{\bN}{\mathbb{N}}
\newcommand{\bR}{\mathbb{R}}
\newcommand{\bQ}{\mathbb{Q}}

\newcommand{\bZ}{\bold{Z}}

\newcommand{\cB}{\mathcal{B}}
\newcommand{\cC}{\mathcal{C}}

\newcommand{\cF}{\mathcal{F}}

\newcommand{\cL}{\mathcal{L}}

\newcommand{\cN}{\mathcal{N}}
\newcommand{\cP}{\mathcal{P}}

\newcommand{\bg}{\bold{g}}

\newcommand{\fC}{\mathfrak{C}}

\newcommand{\Hom}{\text{Hom}}
\newcommand{\Set}{\text{Set}}

\renewcommand{\int}{\text{int}}

\newcommand{\Rmod}{\text{R-Mod}}
\newcommand{\pr}{\text{pr}}
\newcommand{\RMod}{\bK^{\Rmod}}
\newcommand{\Ext}{\text{Ext}}

\newcommand{\comment}[1]{}

\title{A Module-theoretic Introduction to Abstract Elementary Classes}
\author{Will Boney}
\email{wb1011@txstate.edu}
\address{Department of Mathematics, Texas State University, San Marcos, TX, USA}
\date{\today\\Keywords: Abstract Elementary Classes, modules\\The author was supported by the National Science Foundation under grant DMS-2339018.}

\begin{document}

\maketitle
\tableofcontents

\begin{abstract}
The first-order model theory of modules has been studied for decades. More recently, the model theoretic study of nonelementary classes of modules--especially Abstract Elementary Classes of modules--has produced interesting results. This survey aims to discuss these recent results and give an introduction to the framework of Abstract Elementary Classes for module theorists.

\end{abstract}

\section{Introduction}

Model theory is the study of classes of structures that satisfy some fixed set of axioms.  This makes it well-suited to interact with modules, as this is a class of structures that satisfy the module axioms.  Indeed, there has been a rich history of interaction between model theory and module theory, for which \cite{p-modulesbook} is an excellent reference.

Most of this interaction has taken place in the context of first-order model theory.  This is the most well-developed type of model theory, largely due to the Compactness Theorem (Fact \ref{compact-fact}) and other powerful properties of first-order logic.  This is great as several classes of modules are axiomatizable in first-order logic: all modules, torsion-free modules, etc.  To make use of the most powerful tools of model theory, one fixes a \emph{complete} first-order theory of modules (given by specifying invariants conditions \cite[Corollary 2.13]{p-modulesbook}) and looks at all modules over the fixed ring satisfying that theory.  This allows one to make use of the powerful results of classification theory, first laid out by Shelah \cite{sh:c}.

If this was the end of the story, this article wouldn't exist.  On the module-theoretic side, there are many natural classes of modules that are \emph{not} axiomatizable in first-order logic: locally pure injective modules, $\aleph_1$-free abelian groups, flat modules, etc. (see Example \ref{module-ex} for a long list).  On the model theoretic-side, there are many classes of structures that one wants to understand and develop classfication theory for that are not elementary\footnote{A class is \emph{elementary} if it is the models of a complete first-order theory.  A class is nonelementary if it is not elementary, normally only applied to a class of structures in a common language closed under isomorphisms.}.  One could study these classes case-by-case or even by looking at the logics that axiomatize them, but the most popular framework is that of Abstract Elementary Classes (Definition \ref{aec-def}).

Abstract Elementary Classes have an entirely semantic definition, but one that is wide enough to encompass classes axiomatized by a whole host of different logics (see Example \ref{logic-ex}).  Abstract Elementary Classes were introduced by Shelah and first appeared in print in Makowsky \cite[Chapter XX, Definition 1.2.6.(ii)]{bf-model-logic} (where they are called `abstract classes with L\"{o}wenheim number'). Their has been much work developing there classification theory; for references, see \cite{baldwinbook,ramibook, sh:h, bv-survey}.  In the decades since this project began, the focus has mostly been on the abstract development of the theory: fixing an Abstract Elementary Class $\bK$, assuming some additional structural properties (like amalgamation, see Section \ref{basic-prop-ssec}) and potentially some set-theoretic axioms, and developing classification theory.

More recently, there has been an investigation of how this abstract theory interacts with nonelementary classes of modules, especially in the work of Marcos Mazari-Armida and Jan Trlifaj.  This article is designed to introduce the basic ideas and examples in Abstract Elementary Classes, and then describe the applications to modules in the context of the broader theory.

Section \ref{prelim-sec} is dedicated to this first goal.  Recognizing the strange and technical definition of Abstract Elementary Classes, Section \ref{motive-ssec} works to unpack the idea behind what makes a class an Abstract Elementary Classes by working through two basic (and hopefully instructive) examples: modules with submodule, and modules with pure submodule.  Then, a great many examples are given, both from a model-theoretic/logical perspective (Example \ref{logic-ex}) and a module-theoretic perspective (Example \ref{module-ex}).  Finally, we introduce some of the basic notions and properties of Abstract Elementary Classes.

Section \ref{app-sec} is dedicated to the second goal.  It examines three big themes in classification theory: categoricity transfer (Section \ref{cat-ssec}), stability (Section \ref{stab-ssec}) and superstability through the lens of limit models (Section \ref{sstab-ssec}).  For each notion, the abstract theory and motivation is discussed first, and then we outline the applications and interactions of these notions to modules.

We would like to thank John Baldwin, Marcos Mazari-Armida, and the anonymous referee for helpful comments on an initial draft.

\section{Definitions and Examples}\label{prelim-sec}

\subsection{Motivation for definition}\label{motive-ssec}

The axioms for Abstract Elementary Classes (see Definition \ref{aec-def}) fall into the realm of things that are well-motivated and natural once you've internalized them, but are intimidating and off-putting the first time you see them.  In order to give a motivation, we consider two prototypical cases: the class $\RMod$ with the submodule relation $\subset$; and the same class $\RMod$ with the \emph{pure} submodule relation $\subset^{\pr}$.

In each case, we describe a way of analyzing and understanding these cases that will be (mostly) typical of Abstract Elementary Classes.  In the precise language we will provide, both $(\RMod, \subset)$ and $(\RMod, \subset^{\pr})$ are Abstract Elementary Classes with L\"{o}wenheim-Skolem number $\aleph_0+|R|$, while $(\RMod, \subset)$ has the special property of being a \emph{universal class} (see Definition \ref{univ-int-def}.(\ref{univ-cond})).  For those familiar with category theory, the construction we describe is a hands on-version of taking an accessible category $\cK$ with directed colimits, and consider the representation given by the Yoneda embedding $\Hom_{\cK}(\cK_0, -)$, where $\cK_0$ are the generators given by accessibility; this is a concrete way to prove that accessible categories are equivalent to those classes that are axiomatizable in $\bL_{\infty,\infty}$.  Crucially, the classes considered are \emph{finitely} accessible, so this axiomatization is in $\bL_{\infty, \omega}$ and the classes are closed under directed colimits.

The analysis below is not new or surprising, but is designed to highlight the essential features of these classes that mean they satisfy the definition of Abstract Elementary Classes.

\subsubsection{Modules with embeddings} Fix your favorite ring $R$ and consider the category $\RMod$ of $R$-modules with submodule embeddings.  Two warnings before we proceed:
\begin{itemize}
    \item We are concerned with infinite structures.  $R$ might be a finite ring, but each module $M$ we discuss is implicitly assumed to be infinite.
    \item Coming from model theory, the maps between structures are implicitly assumed to be monomorphisms/injections.  This is becasue we want them to preserve the atomic formula `$x\neq y$.'
\end{itemize}
Fix a module $M \in \RMod$.  Given a finite tuple\footnote{We use the model theoretic convention of using bolded variables $\ba, \bb, \bm$, etc. to denote tuples of finite but indeterminate length.} 
$$\ba = a_1, a_2, \dots, a_n \in M$$
we can generate a submodule $\seq{\ba}_M \subset M$.  In this way we can decompose $M$ into a directed union of all of its finitely generated submodules:
$$M \mapsto \{\seq{\ba}_M : \ba \in M\}$$
There are three key points about this decomposition we want to exploit:
\begin{enumerate}
    \item if $\ba \subset \bb$, then $\seq{\ba}_M \subset \seq{\bb}_M$;
    \item (directedness) given $\ba, \bb \in M$, there is $\bc \supset \ba, \bb$; combined with (1), this gives 
    $$M_\ba, M_\bb \subset M_\bc$$
    \item $$M = \bigcup_{\ba \in M} \seq{\ba}_M$$
\end{enumerate}

With these three properties, the correspondence goes the other way!  Fix a directed poset $(I, <)$ and finitely generated $M_i \in \RMod$ for $i\in I$ with the additional property\footnote{This property will be implicitly assumed going forward.} that $i < j$ from $I$ imply that $M_i \subset M_j$.  This allows us to form the directed union of the system
$$M_I := \bigcup_{i \in I} M_i$$
where the universe is the union of the universes.  To define the module operations on $M_I$, we take advantadge of the fact that these operations are finitary (in this case, really binary and unary): given $\ba = a_1, \dots, a_n\in M_I$, there are $i_1, \dots, i_n$ such that $a_{i_k} \in M_{i_k}$.  The directedness of $I$ exactly gives a $i_*\in I$ above each $M_{i_k}$.  Then $\ba \in M_{i_*}$ and any finitary function of $\ba$ can be computed there.

This construction repeats if the models $M_i$ are connected by module embeddings rather than just the submodule relation, that is, a directed system given by 
$$\left(M_i, f_{i,j}\right)_{j<i \in I}$$
such that
\begin{enumerate}
    \item each $f_{i,j}:M_i \to M_j$ is a module embedding; and
    \item the system coheres, which means for all $i<j<k$, the following commutes
\[\begin{tikzcd}
	&& {M_k} \\
	\\
	{M_i} && {M_j}
	\arrow["{f_{i,k}}", from=3-1, to=1-3]
	\arrow["{f_{i,j}}", from=3-1, to=3-3]
	\arrow["{f_{j,k}}"', from=3-3, to=1-3]
\end{tikzcd}\]
\end{enumerate}
While adding a lot of notational complexity, this addition is just cosmetic.

This creates the correspondence in Figure \ref{set-corr-fig}.
\begin{figure}\label{set-corr-fig}\begin{tikzcd}
	\RMod &&& \begin{array}{c} \text{Directed systems of}\\\text{finitely generated}\\ $R$\text{-modules} \end{array}
	\arrow[curve={height=-30pt}, from=1-1, to=1-4]
	\arrow[curve={height=-30pt}, from=1-4, to=1-1]
\end{tikzcd}
\caption{Set correspondence}
\label{set-corr-fig}
\end{figure}
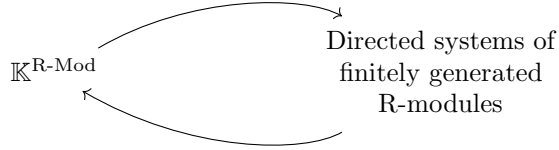

As constructed so far, this correspondence only concerns the modules rather than embeddings between them.  However, we can lift this correspondence to the category level.

Given $f:M\to N$, we can build a system of embeddings between their decompositions: for each $\ba=a_1, \dots, a_k \in M$, define
\begin{eqnarray*}
    \bar{f}_\ba:\langle a_1, \dots, a_n\rangle_M &\to& \langle f(a_1), \dots, f(a_k)\rangle_N\\
    \sum r_i\cdot a_i &\mapsto& \sum r_i \cdot f(a_i)
\end{eqnarray*}
This gives a system $\{\bold{f}_\ba : \ba \in M\}$ such that for all $\ba \subset \bb \in M$, the following commutes
\[\begin{tikzcd}
	{\langle f(\ba)\rangle_N} && {\langle f(\bb)\rangle_M} \\
	\\
	{\langle \ba\rangle_M} && {\langle\bb\rangle_M}
	\arrow[from=1-1, to=1-3]
	\arrow["{\bar{f}_\ba}", from=3-1, to=1-1]
	\arrow[from=3-1, to=3-3]
	\arrow["{\bar{f}_\bb}"', from=3-3, to=1-3]
\end{tikzcd}\]

Putting all of these embeddings together, we get the following commutative diagram

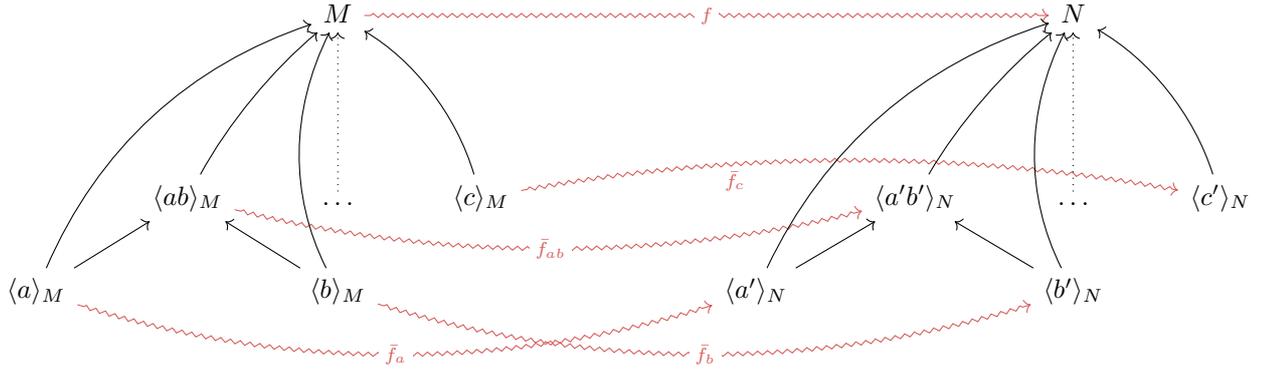
\[\begin{tikzcd}
	&& M &&&&&& N \\
	\\
	\\
	& {\langle ab\rangle_M} & \dots & {\langle c\rangle_M} &&&& {\langle a'b'\rangle_N} & \dots & {\langle c'\rangle_N} \\
	{\langle a\rangle_M} && {\langle b\rangle_M} &&&& {\langle a'\rangle_N} && {\langle b'\rangle_N}
	\arrow["f"{description}, color={rgb,255:red,214;green,92;blue,92}, squiggly, from=1-3, to=1-9]
	\arrow[curve={height=-6pt}, from=4-2, to=1-3]
	\arrow["{\bar{f}_{ab}}"{description}, color={rgb,255:red,214;green,92;blue,92}, curve={height=24pt}, squiggly, from=4-2, to=4-8]
	\arrow[dotted, from=4-3, to=1-3]
	\arrow[curve={height=12pt}, from=4-4, to=1-3]
	\arrow["{\bar{f}_c}"'{pos=0.3}, color={rgb,255:red,214;green,92;blue,92}, curve={height=-18pt}, squiggly, from=4-4, to=4-10]
	\arrow[curve={height=-6pt}, from=4-8, to=1-9]
	\arrow[dotted, from=4-9, to=1-9]
	\arrow[curve={height=12pt}, from=4-10, to=1-9]
	\arrow[curve={height=-24pt}, from=5-1, to=1-3]
	\arrow[from=5-1, to=4-2]
	\arrow["{\bar{f}_a}"{description}, color={rgb,255:red,214;green,92;blue,92}, curve={height=30pt}, squiggly, from=5-1, to=5-7]
	\arrow[curve={height=-18pt}, from=5-3, to=1-3]
	\arrow[from=5-3, to=4-2]
	\arrow["{\bar{f}_b}"{description}, color={rgb,255:red,214;green,92;blue,92}, curve={height=30pt}, squiggly, from=5-3, to=5-9]
	\arrow[curve={height=-24pt}, from=5-7, to=1-9]
	\arrow[from=5-7, to=4-8]
	\arrow[curve={height=-18pt}, from=5-9, to=1-9]
	\arrow[from=5-9, to=4-8]
\end{tikzcd}\]

This means that embeddings between modules can be turned into a coherent system of embeddings between the corresponding directed sequences of models.  This correspondence can also be reversed:  given two directed systems of finitely generated modules $(M_i, f_{i,i'})_{i<i'\in I}$ and $(N_j, g_{j, j'})_{j<j' \in J}$, an embedding between their directed colimits corresponds exactly to the following data:
\begin{enumerate}
	\item an order-preserving map $F_*: I \to J$; and
	\item a coherent system of embeddings $h_i: M_i \to N_{F_*(i)}$, where coherence means the following diagram commutes for all $i<i'\in I$:
\[\begin{tikzcd}
	{M_i'} && {N_{F_*(i')}} \\
	\\
	{M_i} && {N_{F_*(i)}}
	\arrow["{h_{i'}}", from=1-1, to=1-3]
	\arrow["{f_{i,i'}}", from=3-1, to=1-1]
	\arrow["{h_i}"', from=3-1, to=3-3]
	\arrow["{g_{F_*(i),F_*(i')}}"', from=3-3, to=1-3]
\end{tikzcd}\]
\end{enumerate}

This upgrades the correspondence (\ref{set-corr-fig}) between the class of $R$-modules above into a correspondence to the category of $R$-modules with embeddings, as depicted in Figure \ref{cat-corr-fig}.

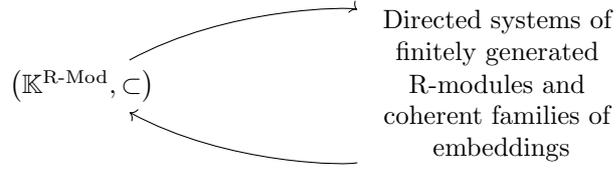
\begin{figure}\label{cat-corr-fig}
\[\begin{tikzcd}
	{\left(\RMod, \subset\right)} &&& \begin{array}{c} \text{Directed systems of}\\\text{finitely generated}\\ $R$\text{-modules and}\\\text{coherent families of}\\\text{ embeddings} \end{array}
	\arrow[curve={height=-30pt}, from=1-1, to=1-4]
	\arrow[curve={height=-30pt}, from=1-4, to=1-1]
\end{tikzcd}\]
\caption{Category correspondence}
\label{cat-corr-fig}
\end{figure}

\subsubsection{Pure embeddings} The context above does not quite provide the generality appropriate for Abstract Elementary Classes.  Given a module $M$ and finite tuple $\ba \in M$, there is a \emph{unique} smallest submodule of $M$ containing $\ba$, namely $\langle \ba \rangle_M$.  However, not every notion of strong subobject that we want to study satisfies this uniqueness property.  Namely, this is not satisfied by the notion of \emph{pure submodule}.  There are many equivalent definitions of pure submodule, but the following is the most useful for our purpose.

\begin{defin}\label{pure-def}
Given $R$-modules $M$ and $N$, we say \emph{$M$ is a pure submodule of $N$}—written $M \subset_\pr N$—iff $M \subset N$ and $M$ reflects the solvability of certain equations from $N$: given a $n\times m$ matrix $A$ of $R$ elements and $\bb \in {}^m M$, 
\begin{eqnarray*}
\text{if $\exists \bx \in {}^n N$ s.t. $A \bx = \bb$,}\\
\text{then $\exists \by \in {}^n M$ s.t.  $A \by = \bb$}
\end{eqnarray*}
\end{defin}

We would like to tell the same story for $(\RMod, \subset_\pr)$ that we did above, but the problem is that there are $\ba \in M$ such that $\langle \ba \rangle_M \not \subset_\pr M$.  For instance, taking $M = \bZ/8\bZ$ and $a = 4$, $M$ thinks that $4$ is divisible by $2$ (expressed as a solution to $`[2][x_1] = [4]’$). But $\langle 4\rangle_M = \{0, 4\}$ does not.  So we cannot find a generated pure submodule covering some parameters.  However, the following means there is \emph{some} pure cover of any given parameter set.

\begin{lemma}\label{pure-dls-lem}
Fix an $R$-module $M$ and let $X \subset M$.  We can build 
$$N \subset_\pr M$$
such that $X \subset N$ and $|N| = |X| + |R| + \aleph_0$.
\end{lemma}

One key part of the lemma is that the bound on the size of $N$ is independent on the size of $N$.  This statement is familiar to logicians as a statement of the Downward L\"{o}wenheim-Skolem Theorem.  The proof of this follows that theorem: build a countable chain of modules $\seq{N_n:n<\omega}$ where $N_0 = \langle X\rangle_M$ and $N_{n+1}$ contains solutions to all matrix equations with parameters in $N_n$; then the union $\bigcup_{n<\omega} N_n$ is the desired $N$.

In Definition \ref{univ-int-def}.(\ref{univ-cond}), a weakening of universality called \emph{closure under intersections} is introduced (but it's mostly what it sounds like).  The class $(\RMod, \subset_\pr)$ is not closed under intersections; thanks to the referee for providing this example.

\begin{example}
    We work with abelian groups (so $R = \bZ$).  Set $N = \bZ(2^\infty)$ to be the\footnote{This example works with any prime $p$.} Pr\"{u}fer 2-group with generators $g_n$ so $2 g_{n+1} = g_n$ and $2g_0 = 0$.

    Define the modules
    \begin{eqnarray*}
    M &=& N \otimes N\\
    M_0 &=& N \otimes \{0\}\\
    M_1 &=& \langle(g_n, 2g_n) : n < \omega\rangle
    \end{eqnarray*}
    Then $M_0, M_1$ are both pure submodules of $M$, as witnessed by the fact that direct summing with $\{0\} \otimes N$ gives $M$.  Then the intersection
    $$M_0 \cap M_1 = \left\{(g_0, 0), (0,0) \right\}$$
    is not a pure submodule of $M$ since $(g_0, 0)$ is not divisible by 2.
\end{example}

Lemma \ref{pure-dls-lem} means that, given a $R$-module $M$, we can decompose it as a directed colimit of small pure submodules, just as before!  There are two crucial changes here:
\begin{enumerate}
    \item We have traded `finitely generated' for `small'. Looking at Lemma \ref{pure-dls-lem}, `small' should mean `of size $|R|+\aleph_0$.'. This cardinal is called the \emph{L\"{o}wenheim-Skolem number} of the class, see Definition \ref{aec-def}.(\ref{ls-cond}).
    \item The new structures are no longer \emph{generated} by anything.  In the submodule case, a map between (free) generators $\ba \in M \mapsto \bb \in N$ lifts uniquely to a map between generated submodules $\langle \ba\rangle_M \to \langle \bb\rangle_N$.  However, no such lift is possible with the pure submodules.
\end{enumerate}

These two conditions make the class $(\RMod, \subset_\pr)$ with pure embeddings much more typical among Abstract Elementary Classes than the class $(\RMod, \subset)$ with submodule embeddings.  Crucially, the following holds:
\begin{itemize}
    \item the category $(\RMod, \subset_\pr)$ is closed under all directed colimits;
    \item the class $\RMod$ satisfies the correspondence Figure \ref{set-corr-fig} with directed systems of small elements of $\RMod$; {\bf however}
    \item the pure embeddings of $(\RMod, \subset_\pr)$ do {\bf not} satisfy the correspondence in Figure \ref{cat-corr-fig} with coherent systems of embeddings between directed systems.
\end{itemize}

This final point exactly follows from the failure of the small models to be generated by finite tuples.  However, we can salvage this result by increasing the directedness of the systems of small models.  Recall that for a cardinal $\kappa$, a partial order $(I, <)$ is \emph{$\kappa$-directed} iff for all $X \subset I$ of size $<\kappa$, there is $i_*\in I$ such that$i < i_*$ for all $i \in X$.  If we upgrade the systems in Figure \ref{cat-corr-fig} to being $(|R|+\aleph_0)^+$-directed, we can recover it.  This construction parallels the $\bK^{up}$ construction in \cite[Lemma 1.24]{sh600}.

\begin{prop}
    Let $R$ be a ring and set $\lambda = \left(|R|+\aleph_0\right)^+$.
    \begin{enumerate}
        \item $M \in \RMod$ corresponds to $\kappa$-directed systems $(M_i, f_{i,i'})_{i<i'\in I}$ where $\|M_i\| = \lambda$ and each $f_{i, i'}$ is a pure embedding.

        Moreover, one such system is given by setting $I = \cP_\lambda (M)$ (the $<\lambda$-sized subsets of $M$) and, for $X \in I$, picking $M_X \subset_\pr M$ to be a pure submodule containing
        $$X \cup \bigcup_{Y\subsetneqq X} M_Y$$
        by Lemma \ref{pure-dls-lem}.
        \item Let $M, N \in \RMod$ correspond to $\kappa$-directed systems $(M_i, f_{i, i'})_{i<i'\in I}$ and $(N_j, f_{j, j'})_{j<j'\in J}$ as above.  Then pure embeddings $f: M \to N$ correspond to order-preserving $F_*:I \to J$ and a coherent system of pure-embeddings $h_i:M_i \to N_{F_*(i)}$.
    \end{enumerate}
\end{prop}

The final condition on an AEC is \emph{coherence} (Definition \ref{aec-def}.(\ref{coh-cond}).  The condition is inherited from model-theoretic contexts where the `strong substructure' contains formulas of great logical complexity.  Then, coherence essentially says that `strong substructure' corresponds to agreeing on the truth of some class of formulas.  When the strongest `strong substructure' relation is purity, this agreement follows for free (see Proposition \ref{free-coh-prop}).

The intuition behind coherence is that strong substructure $M\prec_\bK N$ should represent that, for finite $\ba \in M$, some property holds of $\ba$ in $M$ iff it holds in $N$.  In practice, this condition restricts heavy reference to cardinality in defining AECs.  For instance, attempts to use the equicardinality quantifier in Abstract Elementary Classes normally fails coherence (but see Example \ref{logic-ex}.(\ref{card-quant-cond}) for examples of light references to cardinality that form AECs).

\subsection{Basic definitions}

Now we wish to axiomatize the above behavior into the definition of Abstract Elementary Classes (AECs for short).  AECs were first axiomatized by Shelah \cite{sh88}.  Good references are Baldwin \cite{baldwinbook}, Grossberg \cite{ramibook}, and Boney-Vasey \cite{bv-survey} (the last for tame Abstract Elementary Classes).

Fix a language $\tau$.  We will consider a class $\bK$ of structures along with a strong substructure relationship $\prec_\bK$.  We often write just $\bK$ for the pair $(\bK, \prec_\bK)$.  The following definition captures the above idea that the category $\bK$ (both structures and embeddings) can be captured by $\kappa$-directed colimits for some sufficient $\kappa$; below, $\kappa$ is the cardinal $\LS(\bK)$.

\begin{defin}\label{aec-def}
The pair $(\bK, \prec_\bK)$ is an \emph{Abstract Elementary Classes} iff it satisfies all of the following:
\begin{enumerate}
	\item both $\bK$ and $\prec_\bK$ are closed under isomorphisms;
	\item if $M \prec_\bK N$, then $M \subset_\tau N$;
	\item\label{coh-cond}   (Coherence) given $M_0, M_1, M_2 \in \bK$, if $M_0 \prec_\bK M_2$, $M_1 \prec_\bK M_2$, and $M_0 \subset M_1$, then $M_0 \prec_\bK M_1$;
	\item \label{tv-cond} (Tarski-Vaught axiom) the class is closed under directed colimits which are created in $\Set$;\footnote{The traditional phrasing of this axiom is about closure under unions of chains.  However, a classic result (see \cite[Theorem 2]{ms-poset}) says that these are equivalent conditions.\\
    
    `Created in $\Set$’ means that the underlying set of the colimits is the colimit of the underlying sets.  A crucial contrast is with the category of Hilbert spaces: this category is closed under directed colimits, but the colimit adds more elements than exist in the constituent parts. Lieberman-Rosicky-Vasey \cite[Theorem 18]{lrv-hilbert} go farther and show that there is no way to make Hilbert spaces concrete so directed colimits are created in $\Set$. Metric Abstract Elementary Classes (introduced by Hirvonen and Hyttinen \cite{hh-cathom}) are a formalism appropriate for (Cauchy) complete structures.}  and
	\item\label{ls-cond} (L\"{o}wenheim-Skolem axiom) there is a minimal cardinal $\LS(\bK) \geq |\tau|$ satisfying: given $M \in \bK$ and $X \subset M$, there is $N \in \bK$ such that $X \subset N \prec_\bK$ and $\|N\| = |X| + \LS(\bK)$.
\end{enumerate}
\end{defin}

Some key notation/notions:
\begin{eqnarray*}
	\bK_\lambda &=& \{ M \in \bK : \|M\| = \lambda\}\\
    \bK_{<\lambda} &=& \{ M \in \bK : \|M\| < \lambda\}\\
    \bK_{\geq \lambda} &=& \{ M \in \bK : \|M \| \geq \lambda\}\\
& &( \|M\| \text{ is the size of the underlying set of }M)
\end{eqnarray*}

Note that, due to the Tarski-Vaught axiom, $\bK_\lambda$ uniquely determines $\bK_{\geq \lambda}$ by taking the closure under directed colimits.

A \emph{strong embedding} or \emph{$\bK$-embedding} $f:M \to N$ is a $\tau$-embedding such that$f(M) \prec_\bK N$.  Importantly, {\bf all maps that we write are assumed to be $\bK$-embeddings.}

\subsection{Examples}

As the name implies, this is an abstract concept, so it's nice to give examples.  However, this is perhaps the most difficult section to write.  This difficulty is not  in finding examples (they are numerous!), but in taking the right perspective. It's hard to give a tour of examples that is interesting and meaningful to both logicians and algebraists, so we give two different flavors of examples.  A logician would start by talking about infinitary logic $\bL_{\lambda, \omega}$, explaining what that means, and outlining a few other logics that axiomatize AECs.  I am a logician, so we do that.

In a sense, the algebraic examples are a subset of the logical ones, but they're still useful to list out.

\begin{example}[Logical]\label{logic-ex}\
We begin with a language $\tau$ consisting of some function symbols, relation symbols, and constant symbols.  For module theorists, we typically fix a ring $R$ and use the language\footnote{Including `$-$' and `$0$' is not necessary since they are definable in terms of the other parts of the language, but it simplifies the presentation to include them.} $\tau_{\RMod}$ of $R$-modules:
    \begin{center}
        constant symbol $0$; binary functions $+$ and $-$; unary functions $r\cdot$ for each $r \in R$
    \end{center}
    A logic $\cL$ then takes a language as input and outputs a set of formulas and axioms--more formally called sentences--that a structure can satisfy or not.  When considering AECs, we also want to associate $\cL$ with a strong substructure relation $\prec_\cL^*$ that will turn a class axiomatized by a theory in $\cL$ into an Abstract Elementary Class.

    All\footnote{This is not true of all logics, but is true of all logics considered here and most logics people consider.} logics start with two ingredients:
    \begin{itemize}
        \item \emph{terms} are either free variables, constants in the language, or some composite of the function symbols from the language with inputs as free variables or constants:
        \begin{center}
            $x_1 + 2x_2$ or $1\cdot 0 + y$ 
        \end{center}
        Terms are not true or false statements, but are essentially definitions of functions in the free variables.
        \item \emph{atomic formulas} are relations or equalities between terms:
        \begin{center}
            $x_1 + 2x_2 = 1 \cdot 0 + y$
        \end{center}
        Atomic formulas (like all formulas) are not true or false on their own.  However, we we fix a structure and provide parameters to stand in for any free variables, formulas become true or false.
    \end{itemize}
\begin{enumerate}
    \item \underline{Finitary first-order logic $\bL=\bL_{\omega, \omega}$:} First-order logic--normally written simply as $\bL$ but also written as $\bL_{\omega, \omega}$ when the implicit parameters are being highlighted--is the most basic logic.  We augment atomic formulas by our familiar logical operations:
    \begin{itemize}
        \item negation $\neg$
        \item conjunction $\wedge$ and disjunction $\vee$ (`and' and `or')
        \item existential $\exists$ and universal $\forall$ quantification
    \end{itemize}
    Once a formula has no free variables, it is called an `axiom' or `sentence.' Sentences are either true or false (once a structure has been fixed).

    The familiar sets of axioms for algebraic categories (groups, rings, fields, modules, etc.) are all written in first-order logic.  Stating that a field has characteristic 0 or is algebraically closed is expressible in first-order logic, but requires infinitely many axioms (this is not a problem).

    One of the most important facts about first-order logic is the compactness theorem.

    \begin{fact}[Compactness Theorem, G\"{o}del, Mal'cev]\label{compact-fact}
    Let $T$ be a set of first-order sentences.  Then there is a structure satisfying $T$ iff, for each finite $T_0 \subset T$, there is a structure satisfying $T_0$.
    \end{fact}
    We are focusing on non-first-order contexts where compactness is not relevant, so we will not dwell on this further (except to occasionally notice its absence).
    
    \item \underline{Infinitary logic $\bL_{\lambda, \omega}$:} Fixing an infinite cardinal $\lambda$, $\bL_{\lambda, \omega}$ augments first-order logic by adding infinitary quantification.  Namely, it adds the following logical operation:
    \begin{itemize}
        \item if $\{\phi_i : i < \alpha\}$ is a collection of formulas with $\alpha < \lambda$ and finitely many free-variables among the formulas, then $\bigwedge_{i<\alpha} \phi_i$ and $\bigvee_{i<\alpha} \phi_i$ are also formulas.
    \end{itemize}
    The use of infinitary conjunction and disjunction allows us to assign infinitely much information to a single element.

    For instance, in abelian groups, a formula $\phi(x)$ could only place finitely many restrictions on the order of $x$.  But the logic $\bL_{\omega_1, \omega}$ can include the infinite information that an element is not of finite order in a single formula:
    $$\Phi(x)=\bigwedge_{n < \omega} \left(nx \neq 0\right) $$
    \item \underline{Elementary substructure and Fragments:} In each of the above logics, a natural strong substructure to use is the notion of elementary substructure for the logic, e.g., given $M$ and $N$, 
    \begin{eqnarray*}
        M \prec_{\bL_{\lambda, \omega}} N &\iff& \forall \ba \in M\text{ and } \phi(\bx)\in \bL_{\lambda, \omega}, \text{ we have }\\
        & & M\text{ satisfies }\phi(\ba) \iff N \text{ satisfies } \phi(\ba)
    \end{eqnarray*}

    However, just using this relationship is often much too strong.  For instance, we often study modules with the submodule relation, which is much weaker than the elementary substructure even for first-order logic.

    To remedy this, we instroduce the notion of a \emph{fragment} $\cF$ of a logic, which is simply a collection of formulas that is closed under subformulas; some contexts require that fragments be elementary (closed under the relations of first-order logic), but we use the flexibility of not having this condition.  
    
    For each fragment $\cF$, there is a notion of $\cF$-elementary substructure $M \prec_{\cF} N$ given by only quantifying over the formulas in the fragment.  From this, we can recover:
    \begin{itemize}
        \item the submodule relation as $\prec_{\cF_0}$ for $\cF_0$ the smallest fragment containing all atomic formulas and the module axioms; and
        \item the pure submodule relation as $\prec_{\cF_1}$ for $\cF_1$ the fragment 
        $$\cF_0 \cup \left\{ \exists \bx \phi(\bx, \by) : \phi(\bx, \by) \in \cF_0\right\}$$
    \end{itemize}

    We also use fragments in the context of extra quantifiers, where the fragment substructure relation should also contain whatever extra properties are needed for the strong substructure.
    \item \underline{Extra quantifiers:} Although conjunctions and disjunctions are allowed to be infinitary, quantification is still required to be finitary, that is, we quantify over single elements rather than infinite sequences\footnote{Adding quantification over infinite sequences would run afoul of the Tarski-Vaught axiom (Definition \ref{aec-def}.(\ref{tv-cond})) as new infinite sequences are created in directed colimits.  Adjusting this to allow, for example, quantification over countable sequences and only require the existence of $\aleph_1$-directed colimits leads one to the definition of $\aleph_1$-Abstract Elementary Classes, introduced in \cite{bglrv-muaecs}.}.  However, there are still logical constructions beyond $\bL_{\lambda, \omega}$ that give rise to AECs, primarily in the form of \emph{generalized quantifiers}.  Syntactically, a quantifier takes a formula(s) with free variables and binds one (or more) of those free variables to give a new formula with less free variables.  When viewed this way, our classic existential quantifier turns into the map
    \begin{center}
        `$\phi(x, \by)$' with free variables $x, \by$ $\mapsto$ `$\exists x\phi(x, \by)$' with free variables $\by$
    \end{center}
    This is paired with a semantic rule that inductively defines the satisfaction of the new, quantified formula in terms of the original formula.  For the existential quantifier, this would be written
    \begin{center}
        $M \vDash \text{`}\exists x\phi(x, \bb)\text{'}$ iff there is $a \in M$ such that $M \vDash \text{`}\phi(a, \bb)\text{'}$
    \end{center}
    or
    $$\exists x \phi(x, M) = \left\{\bb : \text{ there is }a \in M\text{ s.t. }(a, \bb) \in \phi(M) \right\}$$
    Each quantifier will come with an associated strong substructure relation that tends to be a strengthening of elementarity\footnote{When I first saw the plethora of strong substructure relations, they first seems to be ad hoc additions to make the AEC axioms hold.  However, two developments speak to the naturality of these additions:
    \begin{itemize}
        \item The most famous example of an `applied' Abstract Elementary Class is Zilber's quasiminimal classes, especially pseudoexponential classes; see \cite{z-pseudoexp}.  These classes can be axiomatized in $\bL_{\omega_1, \omega}(Q_1)$ and the correct natural notion of strong substructure in these examples turned out to be exactly the abstract notion that turns these classes into an AEC.
        \item The extended quantifiers are stronger than the logics $\bL_{\lambda, \omega}$, but are definable in the infinitary quantifier logics $\bL_{\lambda, \kappa}$.  So they can be viewed as including the fragments of this stronger logic that are still tractable.  In each case, the seemingly ad hoc stronger substructure turns out to be exactly the normal elementary substructure definition in the fragment of $\bL_{\lambda, \kappa}$ required to define the quantifier; see \cite[Section 4]{b-cofquant} for more, especially beginning with Example 4.7 there.
    \end{itemize}}:
    \begin{enumerate}
        \item {\bf Cardinality quantifiers:}\label{card-quant-cond} The most basic generalized quantifier is the cardinality quantifiers $Q_\alpha$, which express that certain definable sets are of size $\geq \aleph_\alpha$ (or they are of size $< \aleph_\alpha$ using negation).  Formally, we could express this as  
        \begin{center}
        `$\phi(x, \by)$' with free variables $x, \by$ $\mapsto$ `$Q_\alpha x\phi(x, \by)$' with free variables $\by$
    \end{center}
    where
    $$Q_\alpha x \phi(x, M) = \left\{\bb : \text{ there are at least $\alpha_a$-many }a \in M\text{ s.t. }(a, \bb) \in \phi(M) \right\}$$

    The normal notion of elementary substructure is not enough to make $\bL(Q_\alpha)$-axiomatized class, since we could take an increasing union of models with definable subsets of size $<\aleph_\alpha$ that grow to size equal to $\aleph_\alpha$ in the limit.

    To fix this, when we write `$M \prec_{\bL(Q)}^* N$, we add the condition that, for any negative instance $\neg Q_\alpha x \phi(x, \by)$ in the fragment and $\ba \in M$, we have
    $$\neg \phi(M, \ba) = \neg\phi(N, \ba)$$
    
        \item {\bf Ramsey and equivalent class quantifiers:} There are two quantifiers related to the size of certain sets that are defined in terms of formulas, but are not \emph{definable} in the technical sense and, therefore, can't be captured by the cardinality quantifier:
        \begin{itemize}
            \item The equivalence class quantifier\footnote{This and many other quantifiers are discussed in \cite[Chapters IV and VI]{bf-model-logic}.} $Q^{ec}_\alpha$, where
            $$M\vDash Q^{ec}_\alpha x, y\, \phi(x, y; \ba)$$
            expresses that $\phi(x, y; \ba)$ (with some parameter $\ba$) defines an equivalence relation $x E_\phi y$ with $\geq \aleph_\alpha$-many equivalence classes.
            \item The Magidor-Malitz or (more descriptively) Ramsey quantifiers are $Q^{MM, n}_\alpha$, where
            $$M\vDash Q^{MM, n}_\alpha x_1,\dots, x_n\, \phi(x_1, \dots, x_n; \ba)$$
            expresses that there is a set $X\subset M$ of size $\geq \aleph_\alpha$ such that any $n$-tuple $b_1, \dots, b_n$ from $X$ satisfies $\phi(x_1, \dots, x_n; \ba)$.  Imagining that $\phi$ assigns a coloring the tuple explains the name Ramsey quantifier.
        \end{itemize}
        In both cases, the appropriate strong substructure relation uses the same modification as the cardinality quantifier: if there is a negative instance of the quantifier, then the witnessing set does not grow in the extension.

        Although these notions might not seem immediately connected to module theory, the, e.g., Magidor-Malitz quantifiers have been used in the model theory of modules \cite{b-mm-modules}.
        \item {\bf Cofinality quantifiers:} Cofinality quantifiers were introduced by Shelah \cite{sh43} to find a compact logic beyond finitary first-order.  Cofinality is to linear orders as size/cardinality is to sets: the cofinality of a linear order is the smallest size of a set that is unbounded in the linear order.  For instance, the set of real numbers $\bR$ is uncountable, but the linear order $(\bR, <)$ has countable cofinality since $\bN$ is a countable, cofinal set.

        There are many formalizations of this quantifier, but the simplest is to for a cardinal $\kappa$, define the quantifier $Q^\cf_\kappa$ that takes in two formulas and binds three variables
        $$Q^\cf_\kappa w, x, y \left(\phi(w; \bz), \psi(x, y; \bz)\right)$$
        with the semantics that 
        $$M\vDash Q^\cf_\kappa w, x, y \left(\phi(w; \ba), \psi(x, y; \ba)\right)$$
        iff $\psi(x, y; \ba)$ defines a linear order on the set $\phi(M; \ba)$ with cofinality $\kappa$.  Note that defining a linear order is first-order expressible, so the strength really comes from the statement about the cofinality.    The strong substructure to make Abstract Elementary Classes was given in \cite[Section 3]{b-cofquant}; it is technical so we don't reproduce it here.
        
    \end{enumerate}
\end{enumerate}
\end{example}

But this is not very instructive if you are not steeped in logic.  So we also give several classes of modules that are Abstract Elementary Classes.

{\bf Main idea:} If you’re studying an algebraic property that is \emph{locally verifiable}, the class of structures satisfying that property is probably an Abstract Elementary Classes.

This reduces to another question: what does locally verifiable mean?  One formalism would be to ask if there can be a set of finitary functions that turns the property into a universal property; this can be made precise using Shelah's Presentation Theorem \ref{spt-fact}.
\begin{itemize}
	\item {\bf Locally verifiable:} Locally finite groups\\
	This is a prototypical example: start with a property (finite) and look at all structures where every finitely generated substructure satisfies this property.  These tend to form Abstract Elementary Classes.  
	\item {\bf Not locally verifiable:} Simple groups\\
	Checking if a group $G$ is simple requires searching over it's subgroups to check if they are normal and checking finitely generated subgroups is not enough.  There is a Downward L\"{o}wenheim-Skolem theorem for simple groups that gives a countable simple subgroup of any simple group, however the lack of `finiteness' means that the class is not closed under all unions.
\end{itemize}

Below, we list several natural algebraic classes that form Abstract Elementary Classes.  Unless otherwise stated, these classes form AECs using either the submodule relation or the pure submodule relation.  To streamline presentation, we will list properties of these classes as we introduce them here, even though these properties are not introduced until later (mostly see Section \ref{basic-prop-ssec}).  Most of the examples listed are modules and their variants.

\begin{example}[Algebraic]\label{module-ex}\
Fix a ring $R$.
\begin{enumerate}
    \item \underline{Modules:} The familiar $R$-module axioms are all expressible with just universal quantifiers in finitary first-order logic; call this theory $T_\RMod$. The basic class $(\RMod, \subset)$ of modules with the submodule relation is axiomatized by $T_\RMod$ with the substructure relation.  This fits into the logical framework above by taking $\cF \subset \bL(\tau_\RMod)$ to be the quantifier-free fragment, so $\cF$-elementary substructure is the same as substructure.

    This class is axiomatized in first-order logic and universal, so it satisfies all of the nice properties: amalgamation, joint embedding, no maximal models, $<\omega$-tameness, etc.
    \item \underline{Modules with purity:} Another common class $(\RMod, \subset^\pr)$, which uses pure substructure in place of substructure. Note the underlying models are the same, but we have changed the notion of strong substructure.

    It is not universal (or closed under intersections), but it still satisfies all of the nice properties listed for modules. 
    \item \underline{Torsion-free abelian groups:} The class $(\bK^{tf}, \prec_{pr})$ of torsion-free abelian groups with pure substructure is an Abstract Elementary Class.  Being torsion-free is still first-order axiomatizable, so it satisfies the nice properties listed above \cite[Fact 4.2]{m-desc-limit}
    \item \underline{Torsion modules:} An $R$-module $M$ is torsion if every $x \in M$ is annihilated by some $r \in R$.  When $R$ is infinite this is not first-order, but can be axiomatized in $\bL_{|R|^+, \omega}$ by the sentence
    $$\forall x \bigvee_{r \in R} rx = 0$$
    Although this logic does not satisfy compactness, when $R$ is a PID, a version of an ultraproduct can be developed \cite[Section 5]{b-gamult}.  When we fix a complete first-order theory $T$ and set $\bK^{tor}_T$ to be the models of $T$ that are also torsion modules, this provides enough compactness to prove amalgamation, joint embedding, no maximal models, and that Galois types are syntactic.
     \item \underline{Finitely Butler groups:} Recall that a torsion-free group is \emph{Butler} iff it is a finite rank, pure subgroup of a completely decomposable group.  Then $G$ is \emph{finitely Butler} iff every finitely generated subgroup is a Butler group.  The class $(\bK^{fB}, \subset_{\pr})$ of finitely Butler groups with pure embedding was first considered as an Abstract Elementary Class in \cite[Section 5]{m-desc-limit}.  The proof that this is an AEC is formulaic as a locally verifiable property without unpacking what a Butler group is: if $\cB$ is the class of all finitely generated Butler group with generators, then $\bK^{fB}$ is axiomatized by the $\bL_{|\cB|^+, \omega}$-sentence
    $$\bigwedge_{n<\omega} \forall x_1, \dots, x_n \bigvee_{(G_0; \bg) \in \cB} \bigwedge_{\phi(\bz) \in \qftp(\bg/\emptyset; G_0)} \phi(\bx)$$

    This AEC admits intersections, has joint embedding and no maximal models, and is stable in $\lambda = \lambda^{\aleph_0}$.  Galois types are quanfitier-free syntactic types, so the class is $<\aleph_0$-tame.  It does not have amalgamation, but it has \emph{dense} amalgamation bases (see \cite[Lemma 5.7]{m-desc-limit})\footnotei{Marcos: Does \cite[Lemma 5.7]{m-desc-limit} have a converse: if G is not divisible, does it fail to be an amalgamation base?}.

    \item \underline{Locally finite modules, locally pure injective modules, etc.:} By a similar trick as Finitely Butler groups, the classes of locally finite modules and locally pure injective modules are Abstract Elementary Classes because they are $\bL_{|B|^+, \omega}$-axiomatizable.  Unfortunately, this gives no indication of how to examining the nice properties of the class.  For locally pure injective modules, \cite[Section 3]{m-stab-mod} shows that they satisfy amalgamation, arbitarily large models, and are tame.

    \item\label{root-ext-cond} \underline{Roots of Ext:} Fixing a class $\cN$ of $R$-modules, we can define the corresponding Roots of Ext class:
    \begin{eqnarray*}
        {}^\perp \cN &=& \{M \in \RMod: \Ext^n(M, N) = 0 \text{ for all }N \in \cN\text{ and } 0 < n <\omega\}\\
        M_0 \prec_\perp M_1 &\iff& M_0 \subset M_1 \text{ and } M_1/M_0 \in {}^\perp \cN
    \end{eqnarray*}
    The classes $\left({}^\perp \cN, \prec_\perp\right)$ were first examined as an Abstract Elementary Class by Baldwin, Eklof, and Trlifaj \cite{bet-nperp} and later in \cite{t-aec-tilt, t-roots-ext} and more.  These classes represent an important gap in understanding how to characterize Abstract Elementary Classes: they are one of the very few classes where the strong substrcture relationship is not given syntactically and one of the very few classes where an axiomatization in one of the above logics is not known.

    The class $\left( {}^\perp \cN, \prec_\perp\right)$ is not always an Abstract Elementary Class, but Trlifaj \cite[Theorem 4.4]{t-roots-ext} characterizes the Roots of Ext that are Abstract Elementary Classes as those ${}^\perp \cN$ so ${}^\perp\cN$ is deconstructible (Definition \ref{decon-def}.(\ref{decon-cond}) and closed under direct limits; deconstructibility gives the L\"{o}wenheim-Skolem Axiom \ref{aec-def}.(\ref{ls-cond}) and the direct limits give the colimits propert of the Tarski-Vaught Axiom \ref{aec-def}.(\ref{tv-cond}).

    When Roots of Ext are an Abstract Elementary Class, they satisfy disjoint amalgamation, arbitrarily large models, and have a prime model \cite[Lemma 4.8]{t-roots-ext}.
    
    \item \underline{Flat modules:} The class $(\bK^{flat}, \subset_{\pr})$ of flat modules form an Abstract Elementary Class.  It has amalgamation, joint embedinng, no maximal models, and is tame \cite[Section 6]{lrv-cell}.  If $R$ is a Pr\"{u}fer domain, then this class is closed under intersections.

\end{enumerate}
\end{example}

\subsection{First results}

We list some initial facts about Abstract Elementary Classes before examining deeper properties and results.  The first is a debt from the motivation that the Coherence Axiom is not one that causes problems in ordinary contexts that arise from modules.

\begin{prop}\label{free-coh-prop}
Let $(\bK, \prec)$ be an abstract class such that $\prec$ is a subrelation of $\Sigma_1$-elementarity.  Then $(\bK, \prec)$ satisfy the Coherence axiom, Definition \ref{aec-def}.(\ref{coh-cond}).

In particular, if $\prec$ is pure subgroup or submodule (or some other context that purity makes sense in), then Coherence holds for free.
\end{prop}

Example \ref{logic-ex} gave a great many logics that are guaranteed to axiomatize Abstract Elementary Classes.  A natural question is wether there is a complete logical characterization of which classes are Abstract Elementary Classes.  There is much subtlety to this question (which is discussed in \cite{bv-structural}[Section 4]), but the discussion around Example \ref{module-ex}.(\ref{root-ext-cond}) shows that this kind of exact characterization is not known.  However, from the early days of Abstract Elementary Classes, some connection between syntax and Abstract Elementary Classes has been known.

\begin{fact}[Shelah's Presentation Theorem, {\cite[Lemma 1.8]{sh88}}] \label{spt-fact}
Let $(\bK, \prec_\bK)$ be an Abstract Elementary Class with $\LS(\bK)=\lambda$ in a language $\tau$.  Then there is a an extension $\tau^*\supset\tau$ by $\lambda$-many functions and a  universal $T \subset \bL_{\lambda^+, \omega}(\tau^*)$ theory $T^*$ such that:
\begin{enumerate}
    \item given a $\tau$-structure $M$,
    \begin{eqnarray*}M \in \bK &\iff& \text{ there is an expansion } M^* \text{ of } M \text{ such that }M^*\vDash T^*
    \end{eqnarray*}
    \item given $\tau$-structures $M \subset N$,
    \begin{eqnarray*}
        M\prec_\bK N &\iff& \text{ there are expansions }M^*, N^* \text{ of }M, N\text{ that model } T^* \\
        & &\text{ such that } M^* \subset_{\tau^*} N^*
    \end{eqnarray*}
\end{enumerate}
\end{fact}

Another basic result in Abstract Elementary Classes is the computation of \emph{Hanf numbers} (for model existence).  The question at hand is whether an Abstract Elementary Class must have arbitrarily large models.  In first-order model theory, the Compactness Theorem (\ref{compact-fact}) implies that any elementary class axiomatized by $T$ with a model of size at least $|T|+\aleph_0$ has models of arbitrarily large sizes (this means the Hanf number of such classses is $|T|+\aleph_0$).  Without that theorem, the situation in Abstract Elementary Classes is much more delicate.

Fixing a $\lambda$ to be the $\LS(\bK)$ cardinal, an application of the Axiom of Replacement implies there is a cardinal $H(\lambda)$ such, for any Abstract Elementary Class $\bK$ with $\LS(\bK) \leq \lambda$, if $\bK$ has a model of size at least $H(\lambda)$, then it has arbitrarily large models (Hanf \cite{h-hanf} was the first to make this kind of argument, but not in the context of Abstract Elementary Classes).  This argument shows that such a cardinal exists, but is no help in computing it.  However, for Abstract Elementary Classes, a mostly optimal bound is known:

\begin{fact}\label{hanf-fact}
    Given a cardinal $\lambda$ and an Abstract Elementary Class $\bK$ with $\LS(\bK) \leq \lambda$, if $\bK$ has a model of size $\geq \beth_{(2^\lambda)^+}$ (or even just models of sizes cofinal in this cardinal), then $\bK$ has models of all sizes $\geq \LS(\bK)$.

    Thus, $H(\lambda) \leq \beth_{(2^\lambda)^+}$.
\end{fact}

The proof of this first uses Shelah's Presentation Theorem to show that `the Hanf number of Abstract Elementary Classes with $LS(\bK)\leq \lambda$' is the same as `the Hanf number of classes axiomatized by theories in $\bL_{\lambda^+, \omega}$.'. Then, this second value is known classically by using Morley's Omitting Types Theorem (due to Morley \cite{m-mott} and Chang \cite{c-cpt}, see \cite[Theorem A.3]{baldwinbook} for a proof).

\subsection{Basic Properties} \label{basic-prop-ssec}

The axioms of Abstract Elementary Classes are very bare bones.  On the one hand, this is desirable because this simplicity makes it much easier to determine which classes are AECs.  On the other hand, the minimalist assumptions make it hard (or impossible) to prove deep theorems.  There is a common set of structural properties that have emerged as common hypotheses in theorems, and lemmas to prove when examining specific classes.  As mentioned, we listed the relevant properties earlier when going through examples to streamline the presentation.

The first are a triple of properties that give rise to a \emph{monster model}, a universal domain that we may assume contains all of the modules we want to consider.

\begin{defin} Let $\bK$ be an Abstract Elementary Class.
    \begin{enumerate}
        \item $\bK$ has the \emph{amalgamation property} iff every span has an upper bound, that is, if $M_0, M_1, M_2 \in \bK$ with $M_0 \prec_\bK M_1$ and $M_0\prec_\bK M_2$, then there is $N \in \bK$ and $f_\ell:M_\ell \to N$ for $\ell = 1, 2$ such that the following commutes:
\[\begin{tikzcd}
	{M_1} & N \\
	{M_0} & {M_2}
	\arrow["{f_1}", from=1-1, to=1-2]
	\arrow[from=2-1, to=1-1]
	\arrow[from=2-1, to=2-2]
	\arrow["{f_2}"', from=2-2, to=1-2]
\end{tikzcd}\]
        \item $\bK$ has the \emph{disjoint amalgamation property} iff for every span as above, we can find an upper bound $N$ as there with the additional property that
        $$f_1[M_1] \cap f_2[M_2] = f_1[M_0]$$
        That is, the image of $M_1$ and $M_2$ are disjoint except for the required overlap on the image of $M_0$.
        \item $\bK$ has the \emph{joint embedding property} iff for every $M_0, M_1 \in \bK$, there is $N \in \bK$ and $f_\ell:M_\ell\to N$ for $\ell=0,1$.  In pictures,
\[\begin{tikzcd}
	& N \\
	{M_0} && {M_1}
	\arrow["{f_0}", from=2-1, to=1-2]
	\arrow["{f_1}"', from=2-3, to=1-2]
\end{tikzcd}\]
        \item $\bK$ has \emph{arbitrarily large models} iff for every cardinal $\lambda$, there is $M \in \bK_{\geq \lambda}$.
        \item $\bK$ has \emph{no maximal models} iff for every $M\in \bK$, there is $N \in \bK$ such that $M \precneqq_\bK N$.
    \end{enumerate}
\end{defin}

\begin{remark}
Some comments:
\begin{itemize}
    \item We said three properties and listed four! The final two (`arbitrarily large models' and `no maximal models') are closely related: `no maximal models' implies `arbitarily large models' using the Tarski-Vaught Axiom, Definition \ref{aec-def}.(\ref{tv-cond}) and the converse holds in the presence of the joint embedding property.

    Interestingly (to those looking for more set theory), the converse fails very badly without joint embedding.  Recalling Fact \ref{hanf-fact}, we know an AEC $\bK$ has arbitrarily large models iff it has models of size at least $\beth_{(2^{\LS(\bK)})^+}$.  On the other hand, one can construct an AEC with arbitrarily large models that has no maximal models above a measurable cardinal (if it exists), and Baldwin-Shelah \cite{bs-nmm} have constructed such an example axiomatized by a complete $\bL_{\omega_1,\omega}$ sentence.
    \item Remember each function is a $\bK$-embedding by convention, so $f:M\to N$ always means 
    $$f:M \cong f(M) \prec_\bK N$$
    \item The condition for amalgamation looks very much like a pushout.  While the existence of pushouts implies amalgamation (and are how amalgamation is proved in several classes of modules), amalgamation does not require any uniqueness/colimit properties.
    \item As one might expect, there are many parameterized versions and other variants of these properties ($\lambda$-joint embedding property, density of amalgamation bases,\dots), but we don't list these here.
\end{itemize}
\end{remark}

We can use these properties to build a monster model (Definition \ref{mm-def}).  Before we do this, we introduce an important notion in Abstract Elementary Classes: \emph{Galois types}.  The strength of monster models is best understood by how they simplify this notion.  For those familiar with first-order model theory, a Galois type is a semantic version of the normal syntactic type\footnote{`Type' is best understood here as a synonym of `description.'} that describes the class of elements an element could be mapped to.

\begin{defin}\label{gtp-def}
Fix an Abstract Elementary Class $\bK$.
\begin{enumerate}
    \item A $\bK$-triple $(a, M, N)$ is a collection of $M \prec_\bK N$ with $a \in N$ and\footnote{The condition that $M$ and $N$ are the same size is a restriction to limit the set-theoretic scope of models under consideration; it has no practical effect on the triples we can consider due to the L\"{o}wenheim-Skolem axiom (Definition \ref{aec-def}.(\ref{ls-cond}).} $\|M\|=\|N\|$.
    \item Two triples $(a, M, N)$ and $(b,M, N')$ with the same second coordinate  are \emph{atomically equivalent}--written
    $$(a, M, N) \equiv_{at} (b, M, N')$$
    iff there is an amalgam $f:N\to N''$ and $g:N'\to N''$ with $f(a)=g(b)$ and $f\rest M =g\rest M$.
    \item The \emph{Galois type} of a triple $(a, M, N)$ is written 
    $$\gtp_\bK(a/M; N)$$
    is the equivalence class of $(a,M,N)$ under the transitive closure of atomic equivalance.  The model $M$ is called the \emph{domain} of the type $\gtp_\bK(a/M; N)$.
    \item Given a model $M \in \bK$, the space of Galois types over $M$ is
    $$\S_\bK(M):= \left\{\gtp_\bK(a/M; N) \mid (a, M, N)\text{ is a $\bK$-triple} \right\}$$
    \item Given $M_0 \prec_\bK M$ and $p=\gtp_\bK(a/M; N) \in \S_\bK(M)$, we define the \emph{restriction} of $p$ to $M_0$ to be
    $$p\rest M_0:=\gtp_\bK(a/M_0; N') \in \S_\bK(M_0)$$
    where $N'$ is some/any model $M_0 \prec_\bK N' \prec_\bK N$ such that$a \in N'$ and $\|N'\|=\|N\|$.
\end{enumerate}
In each case, we often omit `$\bK$' when it is clear from context and unambiguous.
\end{defin}

We have implicitly used the basic proposition that atomic equivalence is reflexive and symmetric, so the transitive closure of $\equiv_{at}$ is an equivalence relation.  More importantly, if $\bK$ satisfies amalgamation, then $\equiv_{at}$ is already transitive (see \cite[Exercise 8.8]{baldwinbook}).  The existence of a monster model is often hand-waved using `inessential' large cardinals or other metamathematical tricks (see \cite{hk-monster-model}), but we follow \cite{ramibook} or \cite[Theorem 3.21]{b-taec-notes}  to make this precise inside of $\ZFC$.

\begin{defin}\label{mm-def}
    Fix $\bK$ and a cardinal $\theta \geq \LS(\bK)$.  A \emph{$\theta$-monster model} $\fC$ is a model of $\bK$ satisfying:
    \begin{enumerate}
        \item {\bf $\theta$-universal:} every $M \in \bK_{<\theta}$ embeds into $\fC$;
        \item {\bf $\theta$-model homogeneous:} given $M \prec N$ from $\bK_{<\theta}$ with $f:M \to \fC$, there is $g:N \to \fC$ extending $f$; and
        \item {\bf $\theta$-small types are orbital:} given $M \prec \fC$ of size $<\theta$ and $a, b \in \fC$, we have
        $$\gtp(a/M;\fC) = \gtp(b/M; \fC) \iff \exists f\in \text{Aut}_M \fC\text{ s.t. }f(a) = b$$
    \end{enumerate}
    When $\theta$ is much bigger than the size of models under consideration, we often omit it.
\end{defin}

\begin{fact}
Fix an Abstract Elementary Class $\bK$.  Then the following are equivalent:
\begin{enumerate}
    \item $\bK$ satisifies the amalgamation property, the joint embedding property, and has arbitrarily large models.
    \item For every $\theta$, $\bK$ has a $\theta$-monster model.
\end{enumerate}
\end{fact}

Note that a $\theta$-monster model $\fC$ will be much larger than $\theta$ unless $\theta$ is an inaccessible cardinal.

There are two other relevant properties to define.

One of the key properties of syntactic types in first-order logic follows from the facts that they are sets of formulas and each formula has finitely many parameters.  This immediately gives a local character (even finite character) to syntactic types: any two different syntactic types differ on some finite subset of their domain.  With Galois types having a semantic definition (as orbits in monster models or equivalence classes of amalgams), there is no reason to expect any local character to hold.  However, Grossberg and VanDieren \cite{gv-tamenessone, gv-tamenesstwo, gv-tamenessthree} isolated a local character property of Galois types--called \emph{tameness}--and used it to prove a version of Shelah's Categoricity Conjecture (see Section \ref{scc-sssec}) under the additional hypothesis of a monster model.  Since then, tameness has emerged as a very powerful assumption to make and a property that holds in many algebraic examples; see the survey \cite{bv-survey} for more (and see the ubiquity of tameness in the classes of Example \ref{module-ex}).  As above, there are many parameterizations and variants of this property that we suppress (but are listed in the survey).

\begin{defin} \label{tame-def}
Let $\bK$ be an Abstract Elementary Class.  $\bK$ is \emph{$<\kappa$-tame} iff for every $M \in \bK$ and $p \neq q \in \S_{\bK}(M)$, there is $M_0 \prec M$ with $\|M_0\| < \kappa$ such that
$$p\rest M_0 \neq q\rest M_0$$

We use `$\kappa$-tame' to mean `$<(\kappa^+)$-tame.'

We also allow $\kappa \leq \LS(\bK)$ even thought there might not be any models of size $<\kappa$; in this case, we look for a set of size $<\kappa$ so the restriction of the types are different.\footnote{To make this precise, we would need to give the definition of Galois types over sets, but this is more technical than necessary.  However, we hope the meaning is clear enough without these details.}
\end{defin}

Although Galois types are defined semantically, in nice cases, one can often prove that they are syntactic (see Fact \ref{stable-list-fact} for a partial list).  In this case, the class is $<\omega$-tame for free!

We finish with two properties that are common in algebraic examples: universal classes and closure under intersections.

\begin{defin} \label{univ-int-def}
    Let $\bK$ be an Abstract Elementary Class.  
    \begin{enumerate}
        \item \label{univ-cond} We say $\bK$ is \emph{universal} iff for all $M \in \bK$ and $X \subset M$, the structure $\langle X\rangle_M$ generated by the functions and relations of $M$ is in $\bK$ and satisfies $\langle X\rangle_M \prec_\bK M$.
        \item \label{int-cond} We say $\bK$ \emph{has intersections} iff for all $M \in \bK$ and $X \subset M$, the structure
    $$\cl_M(X):=\bigcap\left\{N \prec_\bK M \mid X \subset N \right\}$$
    is a $\tau$-structure in $\bK$ and satisfies $\cl_M(X) \prec_\bK M$.
    \end{enumerate}
\end{defin}

Recall from the introduction that $(\RMod, \subset)$ is a universal Abstract Elementary Class.  These two properties greatly simplify the conditions necessary to check equality of Galois types whether or not amalgamation holds.

\begin{fact}
Let $\bK$ be an Abstract Elementary class, $M \in \bK$, and $\gtp(a/M; N), \gtp(b/M; N') \in \S_\bK(M)$.
\begin{enumerate} 
    \item If $\bK$ is universal, then
    $$\gtp(a/M; N)=\gtp(b/M; N') \iff \exists f:\langle aM\rangle_N \cong_{M} \langle bM\rangle_N \text{ s.t. }f(a)=b$$
    In this case, $f$ is the unique $\tau$-preserving map generated by 
    $$\id_M \cup \{(a, b)\}$$
    so Galois types are quantifier-free syntactic.
    \item If $\bK$ has intersections, then
    $$\gtp(a/M; N)=\gtp(b/M; N') \iff \exists f:\cl_N(aM) \cong_{M} \cl_{N'}(bM) \text{ s.t. }f(a)=b$$
\end{enumerate}
\end{fact}

Universal classes are exactly those axiomatized by $\bL_{\infty, \omega}$-theories that do not use $\exists$.  There is also a syntactic characterization of Abstract Elementary Classes with intersections, although it is more complex \cite[Corollary 3.11]{bv-structural}.



\section{Modules and Abstract Elementary Classes}\label{app-sec}

\subsection{Categoricity and Deconstruction}
\label{cat-ssec}
\subsubsection{Categoricity transfer in abstract Abstract Elementary Classes} \label{scc-sssec} One of the biggest test questions in Abstract Elementary Classes is \emph{Shelah’s Categoricity Conjecture} and its many variants.  

The roots of this conjecture go back to a classic result in linear algebra: vector spaces have a nice invariant—the dimension—such that any two vector spaces (over the same field) are isomorphic iff they have the same dimension. For “small” dimensions, the dimension of the vector field can be “hidden” by its cardinality, but for large dimensions, the cardinality and the dimension of the vector space are the same.  This leads to the following observation:
\begin{center}
Given any two (real) vector spaces, if they are of same size larger than $|\mathbb{R}|$, then they are isomorphic.
\end{center}
The behavior of ‘any two structures of the same size being isomorphic’ is called \emph{categoricity in power}\footnote{`Power’ is an old-fashioned term for size or cardinality.  Also, the term `categoricity’ once meant that \emph{every} structure in a class was isomorphic, in the sense that all of the models were required to be isomorphic regardless of size.  However, in the context of model theory, this is very uncommon, so `categoricity’ with no modification is normally used to mean `categoricity in power.’}.  

In 1954, \L o\'{s} observed that, among countable first-order theories, there tended to be an all or nothing behavior with categoricity in power, and conjectured for the class of models of a countable first order theory, either the theory was categorical in \emph{every} uncountable cardinal or it was categorical in \emph{no} uncountable cardinals; the question of countable categoricity turns out to be independent from the rest.

Morley \cite{m-categoricity} proved this in the affirmative in his 1964 PhD thesis.  This was followed by work of Baldwin and Lachlan \cite{bl-categoricity} that showed that this result can be explained by finding—in a class of models of a countable theory that is categorical in an uncountable cardinal—admit a dimension invariant that determines the isomorphism class of the model.  Shelah \cite{sh31} further extended Morley’s result showing that it also holds amongst \emph{uncountable} theories, although the cutoff must be larger.  

These results can be summarized in the following result (phrased in a way suggestive of future generalizations):

\begin{fact}[{Morley \cite{m-categoricity}, Baldwin-Lachlan \cite{bl-categoricity}, Shelah \cite{sh31}}]\
For every $\lambda$, there is some $\mu_\lambda$ such that for any first-order theory of size $\lambda$, if the models of $T$ are categorical in \emph{some} $\kappa \geq \mu_\lambda$, then the models of $T$ are categorical in \emph{every} $\kappa' \geq \mu_\lambda$.

In particular, we can take $\mu_\lambda =\lambda^+$.  If $\lambda=\omega$, then there is a dimension associated with the models of $T$ that determines the isomorphism class (including for countable models).\footnote{Unfortunately, Shelah’s extension to uncountable theories does {\bf not} come with an associated dimension.  This means that the motivating example from real vector spaces is not actually explainable by the abstract model theory.  However, the linear dimension for vector spaces of countable field (such as $\bQ$) or the transcendence dimension for algebraically closed fields of a fixed characteristic is discoverable via the Baldwin-Lachlan result).}
\end{fact}

Morley’s theorem and the work that followed it lead to a renaissance in first-order model theory, especially Shelah’s classification theory, first laid out in \cite{sh:a}.  Highlights of the successes of this program include Hruskovski's proof of the function-field Mordell-Lang Conjecture \cite{h-mordell-lang}, improvements to Szemer\'{e}di's Regularity Lemma \cite{ms-szemeredi}, work on definable groups \cite{hp-pseudo}, and contributions to differential algebra \cite{fs-j}.

Shelah’s Categoricity Conjecture is the generalization of Morley’s Categoricity Theorem to non elementary classes.  It was initially stated for theories in $\bL_{\omega_1, \omega}$ or $\bL_{\lambda, \omega}$, but is now normally studied in the context of Abstract Elementary Classes.

\begin{conj}[Shelah’s Eventual Categoricity Conjecture]
For every $\lambda$, there is some $\mu_\lambda$ such that for any Abstract Elementary Class $\bK$ with $\LS(\bK) = \lambda$, if $\bK$ is categorical in \emph{some} $\kappa \geq \mu_\lambda$, then $\bK$ is categorical in \emph{every} $\kappa'\geq \mu_\lambda$.
\end{conj}

The standard version (with no `Eventual') speicifies that $\mu_\lambda = \beth_{\left(2^\lambda\right)^+}$ (which is the Hanf number $H(\lambda)$), but even proving the more general eventual version stated above would be a huge advance.

To give context for the results in AECs of modules, we list some results for abstract AECs:

\begin{fact}\label{scc-fact}
    Let $\bK$ be an Abstract Elementary Class with $\LS(\bK) = \lambda$.
    \begin{enumerate}
        \item Suppose $\bK$ has amalgamation.  Set $\mu_\lambda:= \beth_{\left(2^{\beth_{\left(2^{\lambda}\right)^+}}\right)^+}$.  If $\bK$ is categorical in some $\kappa^+ >\mu_\lambda$, then $\bK$ is categorical in all $\kappa'$ such that $\mu_\lambda \leq \kappa' < \kappa^+$. (\cite[Theorem 2.7]{sh394}, Downward Categoricity Transfer from a Successor)
        \item \label{scc-frame-cond} Suppose that $2^{\kappa^{+n}} < 2^{\kappa^{+(n+1)}}$ for all $n<\omega$.  If $\bK$ is categorical in $\kappa^{+n}$ for all $n<\omega$, then $\bK$ is categorical in all $\kappa' \geq \kappa$ (\cite[Conclusion 12.43]{sh705})
        \item \label{scc-tame-cond} Suppose $\bK$ has amalgamation, no maximal models and is $\kappa$-tame.  Set $\mu_\lambda = \lambda^++\kappa$. If $\bK$ is categorical in $\kappa^+ > \mu_\lambda$, then $\bK$ is categorical in all $\kappa ' > \kappa^+$. (\cite[Theorem 5.2]{gv-tamenessthree}, Upward Categoricity from a Successor with Tameness)
        \item Suppose there is a strongly compact cardinal $\mu$ above $\lambda$.  Set $\mu_\lambda = \mu$.  If $\bK$ is categorical in some $\kappa^+ >\mu_\lambda$, then $\bK$ is categorical in all $\kappa' \geq \mu_\lambda$. (\cite[Theorem 7.4]{b-tamelc}, Shelah's Eventual Categoricity Conjecture for Successors from a Strongly Compact)
        \item \label{scc-sv-cond} Suppose that $2^\lambda < 2^{\lambda^+}$ for all $\lambda$ and $\bK$ has amalgamation and arbitrarily large models.  Set\footnote{This result is phrased this way for consistency with other results, the result is actually much more powerful.} $\mu_\lambda = \beth_{\left(2^\lambda\right)^+}$.  If $\bK$ is categorical in some $\kappa \geq \mu_\lambda$, then $\bK$ is categorical in all $\kappa' \geq \mu_\lambda$. (\cite[Corollary 9.7]{v-catspec} using results from \cite{sv-multidim}, see also \cite[Corollary 10.14]{v-catspec})
    \end{enumerate}
\end{fact}

Note that additional assumptions are typically needed, either on the AEC or on the set-theoretic universe.  In these partial results, note that the threshold cardinal is quite high.

We should also comment on the description of Shelah's Categoricity Conjecture as a `test question.'. This is used to indicate that it is not expected that a positive solution will lead to new results by citing the true statement.  Instead, proving the conjecture is a \emph{test} of the robustness of the structure theory of Abstract Elementary Classes, and it is the \emph{techniques} used to solve it that will be used to prove other theorems.  This has been seen in the partial results listed above, e.g., tameness in Fact \ref{scc-fact}.(\ref{scc-tame-cond}) or good $\lambda$-frames in Fact \ref{scc-fact}.(\ref{scc-frame-cond}). 

\subsubsection{Categoricity in Abstract Elementary Classes of Modules}

We now situate categoricity in the study of Abstract Elementary Classes of modules.  An important point is that categoricity is solely a property of the models in an Abstract Elementary Class $\bK$ and is \emph{not} impacted by the particular choice of strong substructure $\prec_\bK$ (although some strong substructure must be given to witness that the collection of modules form an Abstract Elementary Classes).  One theme throughout the results is that the threshold cardinals (which we normally denote $\mu_\lambda$) are \emph{much} lower in the module examples than in general: compare the bound $|R|+\aleph_0$ in the clauses of Fact \ref{scc-list-fact} below with the corresponding bound $\beth_{(|R|+\aleph_0)^+}$ in Fact \ref{scc-fact}.(\ref{scc-sv-cond}) above.

Categoricity in AECs of modules have been studied by Mazari-Armida \cite{m-cat}, Saroch and Trlifaj \cite{st-decon}, and Trlifaj \cite{t-cat-ext}.  In each of these contexts, the hypothesis and the details change, but the argument tends to follow a similar roadmap:
\begin{enumerate}
    \item Identify some type of special modules 
    \item Use categoricity in some $\lambda$ to show every model of size $\lambda$ must be isomorphic to some direct sum of a special model (or some other kind of structure result).
    \item This will imply there is a unique special model $M$ and that every model in $\bK$ is isomorphic to a direct sum of $M$
\end{enumerate}

Mazari-Armida sets up a framework for this when special means \emph{strongly indecomposable}.  An $R$-module is strongly indecomposable if its endomorphism ring is local.  Crucially, the Krull-Schmidt-Remak-Azumaya Theorem (see \cite[2.1]{f-module}) provides that a decomposition of a module $M$ into strongly indecomposables is unique up to a permutation of the indices; however, this does not provide that a decomposition exists.  Thus we call a class of $R$-modules $\bK$ \emph{nice} iff every module can be written as a direct sum of strongly indecomposables and there is an upper bound $\mu_\bK$ on the size of the strongly indecomposables (this is \cite[Hypothesis 2.4]{m-cat}).

Given a module $M$ and a cardinal $\lambda$, we use $M^{(\lambda)}$ to denote the direct sum $\oplus_\lambda M$ of $\lambda$-many copies of $M$.

\begin{prop}
    If $\bK$ is a nice class of modules, then $(\bK, \subset)$ is an AEC with L\"{o}wenheim-Skolem number $|R|+\mu_\bK+\aleph_0$.
\end{prop}


In nice classes, there is a tight relation between categoricity and the number $\mu_\bK$ of strongly indecomposable modules.

\begin{lemma}[{\cite[Lemma 2.7]{m-cat}}]\label{cat-nice-lem}
Let $\bK$ be a nice class of modules and $\kappa > \mu_\bK$.  Then
\begin{center}
    $\bK$ is categorical in $\kappa$ iff $\bK$ has exactly one strongly indecomposable module (up to isomorphism)
\end{center}
\end{lemma}

{\bf Proof:} First, suppose $\bK$ has exactly one strongly indecomposable $M$.  The direct sum $M^{(\mu)}$ has size $\mu+ \|M\|$, so every $N \in \bK_\kappa$ is isomorphic to $M^{(\kappa)}$.

Second, suppose $\bK$ is categorical in $\kappa > \mu_\bK$.  Let $M, N$ be strongly indecomposable.  Then $M^{(\kappa)}$ and $N^{(\kappa)}$ are of size $\kappa$ and therefore isomorphic by the categoricity assumption.  The unique decomposability of the Krull-Schmidt-Remak-Azumaya Theorem implies that $M$ and $N$ are isomorphic.\hfill \dag\\

Since the right hand implication of Lemma 2.7 does not depend on the cardinal $\kappa$, this shows that Shelah's Categoricity Conjecture holds for nice classes of modules (with the much lower threshold cardinal).

There are many examples of classes of modules that turn out to be nice.  In many cases, the categorical side of this dichotomoy (where $\bK$ is categorical in some/all $\lambda > \mu_\bK$) also comes with a characterization of the ring $R$.  The various definitions of the ring properties are not crucial for our purposes, but are available at the citations.

\begin{fact}\label{scc-list-fact}
Fix an associative ring $R$ with unity.  Each of the following classes is nice and, therefore, satisfies Shelah's Categoricity Conjecture.
\begin{enumerate}
    \item $\RMod$\\
    This class is categorical in some power iff $R$ is the ring of matrices over some division ring. \cite[Lemma 2.12]{m-cat}
    \item the class of locally pure-injective modules with pure embeddings\\
    This class is categorical in some power iff $R$ is a full matric ring over a division ring (and thus reduces to the previous case). \cite[Theorem 3.4]{m-cat}
    \item \label{abs-pure-case} the class of absolutely pure modules\\
    This class is categorical in some power iff $R$ is left Noetherian and there is a unique indecomposable injective module up to isomorphism. \cite[Theorem 3.11]{m-cat}
    \item the class of locally injective modules\\
    This class is categorical in some power iff $R$ is left Noetherian and there is a unique indecomposable injective module up to isomorphism. \cite[Theorem 3.19]{m-cat}
    \item the class of flat modules\\
    This class is categorical in some power iff $R$ is left perfect and $R/J(R)$ is the ring of matrices over some division ring iff $R$ is the ring of matrices over a local ring with left $T$-nilpotent maximal ideal. \cite[Theorem 3.26]{m-cat}
\end{enumerate}
\end{fact}

The extra information about $R$ in these cases comes from showing that relevant properties interact nicely with direct sums and any special module we care about is the direct sum of our strongly indecomposable generator.  We give an example below. Recall that if every absolutely pure left $R$-module is injective, then $R$ is left Noetherian.

{\bf Proof of Fact \ref{scc-list-fact}.(\ref{abs-pure-case}):}\footnotei{To Marcos: In, e.g., \cite[Theorem 3.11]{m-cat}, you use some high powered Vasey result to get categoricity in certain strong limits, and make the transfer there.  However, it looks like you independently show categoricity everywhere.  Do you really need the Vasey result?}Suppose this class is categorical in some power.  Then it is categorical in every power by Lemma \ref{cat-nice-lem}.  Let $M$ be an absolutely pure $R$-module.  Then its injective envelope $E(M)$ is absolutely pure.  By categoricity,
$$E(M) \cong M^{(\mu)}$$
for $\mu$ being the size of $E(M)$.  That means $M^{(\mu)}$ is injective and direct sums preserve non-injectivity, so $M$ is injective.\hfill \dag\\

Note that many of the classes in Fact \ref{scc-list-fact} satisfy nice AEC properties like amalgamation, tameness, etc. \cite[Lemma 3.5 and 3.10]{m-stab-mod}

Another direction \cite{st-decon, t-cat-ext} for categoricity focuses on \emph{strong splitters} (Definition \ref{decon-def}.(\ref{ss-cond}))as the special models in the context of \emph{deconstructible classes} of modules.

\begin{defin}\label{decon-def}
Fix a ring $R$, an infinite cardinal $\kappa$, and a collection $\cC$ of $\kappa$-presentable modules.
\begin{enumerate}
    \item Given an module $M$, a \emph{$\cC$-filtration of $M$} is a continuous, increasing chain $\{M_i: i < \alpha\}$ where
    \begin{itemize}
        \item $M_0 = 0$;
        \item $M = \bigcup_{i<\alpha} M_i$; and
        \item for each $i<\alpha$, $M_{i+1}/M_i$ is isomorphic to an element of $\cC$.
    \end{itemize}
    \item \label{decon-cond} The class $\bK$ is \emph{deconstructible} iff there is some $\kappa$ such that $\bK$ is precisely the collection of all modules with a $\bK_{*<\kappa}$-filtration, where $\bK_{*<\kappa}$ are the $<\kappa$-presented modules in $\bK$.  In this case, we say $\bK$ is \emph{generated} by $\bK_{*<\kappa}$.
    \item \label{ss-cond} An $R$-module $M$ is a \emph{strong splitter} iff
    $$\Ext_R^1(M,M^{(\kappa)}) = 0$$
    for some $\kappa\geq|R|+\aleph_0$ such that $M$ is $\leq\kappa$-presented.\footnote{Note the definition requires this condition for \emph{all} direct sums, but \cite[Lemma 2.2]{t-cat-ext} shows that this single set is enough.}.
\end{enumerate}

Note that $M$ being a strong splitter means that given any $N \supset M^{(\kappa)}$ such that $N/\left(M^{(\kappa)}\right) \cong M$, we must have a submodule $Q \subset N$ such that $M^{(\kappa)} \oplus Q = N$, that is, $M^{(\kappa)}$ splits in $N$.
    
\end{defin}

Importantly, not all deconstructible classes of modules are Abstract Elementary Classes (implicitly using the submodule relation unless we state otherwise) and the following results cover \emph{all} deconstructible classes.  However, there are many examples of deconstructible classes that are Abstract Elementary Classes: any Roots of Ext (Example \ref{module-ex}.(\ref{root-ext-cond})) that form an AEC are deconstructible \cite[Theorem 2.2]{st-decon}.

In deconstructible classes, strong splitters are closely tied to categoricity.  In fact, the existence of a single non-strong splitter massively restricts the categoricity spectrum. A helpful observation is that it suffices to check a particular model to see if there are any non-strong splitters.

\begin{fact}[{\cite[p. 375]{t-cat-ext}}]
Fix a deconstructible class $\bK$ that is generated by $\cC$.  Define
$$M_\cC = \bigoplus_{N \in \cC} N$$
If $M_{\cC}$ is a strong splitter, then every $M \in \bK$ is a strong splitter.
\end{fact}

A single model (and therefore $M_\cC$) failing to be a strong splitter already means that there could be at most one categoricity cardinal!

\begin{fact}[{\cite[Lemma 2.8 and Theorem 2.12]{t-cat-ext}}]
Let $\bK$ be a deconstructible class of $R$-modules.  Suppose there is a non-strong splitter $M \in \bK$ of size $\kappa \geq |R|+\aleph_0$.  Then $\bK$ is categorical in at most one cardinal $\lambda \geq \kappa$ and such a $\lambda$ is either equal to $\kappa$ or singular.
\end{fact}

Of course, these restrictions on $\lambda$ beg the question if such a $\lambda$ could possibly exist.  If $\bK$ is $\mu^+$-deconstructible with $\mu\geq |R|+\aleph_0$, then such a $\lambda$ must be $\leq 2^\mu$ by \cite[Theorem 4.1]{st-decon}.  This question is open \cite[p. 381]{t-cat-ext}.

{\bf Proof Sketch:} The key in both cases is to use the non-strong splitting of some $M$ and potentially categoricity in $\lambda$ to build models of size $M^*_\mu$ of size $\mu$ that \emph{cannot} be isomorphic to $M^{(\mu)}$, which is another model of size $\mu$ in $\bK$.

To do this, build a continuous, increasing chain $\{M_\alpha: \alpha < \mu\}$ starting with $M_0 = M^{(\lambda)}$ and using the following successor step: if $M_\alpha$ is isomorphic to some $\oplus_\chi M$ (necessarily with $\chi \geq\lambda$), then the definition of not strongly splitting means that there is an extension $M_{\alpha+1} \supset M_\alpha$ with $M_{\alpha+1}/M_\alpha \cong M$ and $M_\alpha$ does not split in $M_{\alpha+1}$.  The first condition ensures $M_{\alpha+1}$ (and the eventual union) is in the deconstructible class, and the second condition will force us away from the union being a direct sum itself.

Set $M^*_\mu = \bigcup_{\alpha < \mu} M_\alpha$.  The proof that $M^*_\mu$ is not isomorphic to $\oplus_\mu M$ is different in each case, but the essential argument is to a back and forth argument from a supposed bijection that contradicts our construction: it finds a stage $M_\beta \subset M_{\beta+1}$ such that $M_\beta$ is a direct summand of $M$ and, hence, $M_\beta$ splits in $M_{\beta+1}$.\hfill \dag\\

On the structure side, we may assume that all models are strong splitters.  Then categoricity reduces to finding a power with all direct sums of generating modules to that power are isomorphic.

\begin{fact}[{\cite[Lemma 2.9]{t-cat-ext}, \cite[Lemmas 4.3 and 4.5]{st-decon}}]
Let $\bK$ be deconstructible and generated by $\cC$ such that $M_\cC$ is a strong splitter.  Set
$$\nu = \aleph_0+|R| + |M_\cC|$$
Then for each $\lambda > \nu$
\begin{center}
    $\bK$ is categorical in $\lambda$ iff every $M, N \in \cC$ satisfies $M^{(\nu)} \cong N^{(\nu)}$.
\end{center}
    
\end{fact}

Again, the right hand condition contains no $\lambda$, so this gives an instance of Shelah's Categoricity Conjecture.

{\bf Proof Sketch:} First, suppose that $\bK$ is categorical in $\lambda$.  Then $M^{(\lambda)}\cong N^{(\lambda)}$, and we want to bring this down to $\nu$.  Since $M^{(\nu)}$ is a direct summand of $M^{(\lambda)}$, it is one of $N^{(\lambda)}$, so $M^{(\nu)} \oplus M_*\cong N^{(\nu)}$ for some module $M_*$.  Then we use a version of Eilenberg's trick to show
\begin{eqnarray*}
    N^{(\nu)} \cong \left(N^{(\nu)}\right)^{(\omega)} \cong \left(M^{(\nu)} \oplus M_*\right)^{(\omega)} \cong M^{(\nu)} \oplus \left(M_* \oplus M^{(\nu)}\right)^{(\omega)} \cong M^{(\nu)} \oplus N^{(\nu)}
\end{eqnarray*}
Dually, $N^{(\nu)}$ is a direct summand of $M^{(\nu)}$, so we can show
$$M^{(\nu)} \cong M^{(\nu)} \oplus N^{(\nu)}$$
Thus, $M^{(\nu)}\cong N^{(\nu)}$.

Second, assume the right side and take some $N_0 \in \cC$ with $\|N_0\| \leq \nu < \lambda$.  Fix $M \in \bK_\lambda$ and we will show this is isomorphic to $N_0^{(\lambda)}$.  Since every module in $\cC$ is a strong splitter, we have that
$$M \cong \bigoplus_{N \in \cC} N^{(\kappa_N)}$$
for some cardinals $\kappa_N \leq \lambda$.  We partition $\cC$ into
\begin{eqnarray*}
    \cC_0 &=& \left\{ N \in \cC : \kappa_N < \nu\right\}\\
    \cC_1 &=& \left\{ N \in \cC : \kappa_N \geq \nu \right\}
\end{eqnarray*}
and set
\begin{eqnarray*}
    M_0 &:=& \bigoplus_{N \in \cC_0} N^{(\kappa_N)}\\
    M_1 &:=& \bigoplus_{N \in \cC_1} N^{(\kappa_N)}\\
    M &\cong& M_0 \oplus M_1
\end{eqnarray*}
For $M_0$, we have $|\cC_0|\leq |\cC| < \lambda$, so $\|M_0\| < \lambda$ and $\|M_1\|=\lambda$.  For $M_1$, we use our assumption and the fact $\kappa_N \geq \nu$ to conclude
\begin{eqnarray*}
    M_1 &=& \bigoplus_{N \in \cC_1} N^{(\kappa_N)} = \bigoplus_{N \in \cC_1} \left(N^{(\nu)}\right)^{(\kappa_N)}\\
    &=&\bigoplus_{N \in \cC_1} \left(N_0^{(\nu)}\right)^{(\kappa_N)} = \bigoplus_{N \in \cC_1} N_0^{(\kappa_N)}\\
    &=& \bigoplus_{\sum_{N\in \cC_1}\kappa_N} N_0 = \bigoplus_\lambda N_0\\
    &=& N_0^{(\lambda)} 
\end{eqnarray*}
By our assumption, $M_0$ is isomorphic to a direct summand of $N_0^{(\lambda)}$ and $M \cong M_0 \oplus N_0^{(\lambda)}$.  So Eilenberg's trick can be adapted (see the proof of \cite[Lemma 4.5]{st-decon}) to show
$$M \cong N_0^{(\lambda)}$$\hfill \dag\\

\subsection{Galois Stability in Modules}\label{stab-ssec}

\subsubsection{Stability in Model Theory} One of the basic dividing lines of classification is stability theory.  In first-order model theory, stability has several equivalent characterizations.

\begin{fact}\label{fo-stab-fact}
Let $T$ be a complete first-order theory.  The following are equivalent:
\begin{enumerate}
    \item \label{fo-syn-cond} {\bf Syntactic stability:} There is some $\lambda$ such that, for any $M \vDash T$ of size $\lambda$, $|S(M)| = \lambda$, where $S(M)$ denotes the set of \emph{syntactic} types over $M$.
    \item {\bf Undefinability of order:} There is no formula $\phi(\bx, \by)$ that defines an infinite linear order on some elements from a model of $T$.
    \item {\bf Stable nonforking:} There is an independence relation satisfying certain properties (see, for instance, \cite[Theorem 8.5.10]{tz-model-theory}).
\end{enumerate}
\end{fact}

Stability is a powerful tool in model theory and classification theory. Relevant to modules, every first-order theory of modules is stable \cite{b-cat-mod, f-thesis}, and Prest's book \cite{p-modulesbook} is a good introduction to the first-order model theory of modules.

When moving the discussion to Abstract Elementary Classes, the definitions shift to semantic versions.  The type counting version of stability (Fact \ref{fo-stab-fact}.(\ref{fo-syn-cond})) is normally taken as the definition of stability in first-order, and we mimic this in AECs.

\begin{defin}
    Let $\bK$ be an Abstract Elementary Class.  We say $\bK$ is \emph{$\lambda$-Galois stable} iff for every $M \in \bK_\lambda$, 
    $$|\S_\bK(M)| = \lambda$$
    We say $\bK$ is \emph{Galois stable} iff it is $\lambda$-Galois stable for some $\lambda$.
\end{defin}

Note that, under amalgamation, an equivalent definition of $\lambda$-Galois stability is the existence of $(\lambda, \alpha)$-limit models; see Section \ref{sstab-ssec} and Proposition \ref{exist-limit-prop} for more.

\subsubsection{Galois stability in modules} The most basic question about stability in AECs of modules is whether the analogue of the Baur-Fisher result holds for Galois stability, first asked by Mazari-Armida \cite[Question 2.12]{m-fuchs}.
\begin{question} \label{marco-quest}
    Let $R$ be an associative ring with unity and $\bK \subset \RMod$ such that $(\bK,\subset_{\pr})$ is an Abstract Elementary Class.  Is $(\bK, \subset_{\pr})$ Galois stable?
\end{question}

All evidence points to a positive answer.  What follows is a long list of examples of classes of modules with pure submodule that have been proven to be Galois stable.  Of course, the fact this question is open indicates there is no known example that is not Galois stable.

\begin{fact}\label{stable-list-fact}
    Let $R$ be an associative ring with unity.  The folllowing classes $\bK$ of modules give that the Abstract Elementary Class $(\bK, \subset_{\pr})$ that is stable.
    \begin{enumerate}
        \item \label{list-mod-cond} $\RMod$ (follows from first-order since Galois types are syntactic)
        \item \label{list-tf-cond} $\bK^{tor}$, the class of torsion abelian groups (\cite[Corollary 4.8]{m-fuchs})
        \item the torsion modules satisfying some complete theory of modules {\bf when} $R$ is a PID (\cite[Theorem 5.17]{b-gamult}, Galois types are syntactic)
        \item reduced torsion-free Abelian groups (\cite[Claim 1.2]{sh820})
        \item flat modules (\cite[Theorem 4.3]{lrv-cell}, uses stable independence relation)
        \item definable classes of modules (\cite[Theorem 3.16]{km-univ}, Galois types are syntactic)
        \item $\bK$ closed under direct sums and pure-injective envelopes (\cite[Theorem 3.11]{m-stab-mod}, Galois types are syntactic)
        \item $\bK$ closed under direct sums, pure submodules, and pure epimorphic images (\cite[Theorem 4.17]{m-stab-mod}, uses stable independence relation)
        \item $\bK$ is closed under submodules and has arbitrarily large modules {\bf when} $R$ is von Neumann regular (\cite[Lemma 5.10]{m-stab-mod}, Galois types are syntactic)
    \end{enumerate}
\end{fact}

As listed, there are two common methods of proof.  First, if one can prove that Galois types are syntactic for some fragment of first-order logic (typically either quantifier-free of p.p. formulas), then the number of Galois types is bounded above by the number of first-order syntactic types, so Galois stability follows from the classic Baur-Fischer result for first-order. Second, some of the properties of a stable independence relation can be used to prove Galois stability using a counting argument.

Note that the first-order Baur-Fisher result deals with the number of syntactic types in a \emph{complete} first-order theory of modules.  While some of the examples are refinements of these elementary classes (Fact \ref{stable-list-fact}.(\ref{list-tf-cond}) for instance), others are not (Fact \ref{stable-list-fact}.(\ref{list-mod-cond}) for instance).  There is some subtle interplay of monotonicity when it comes to Galois types, but it is not the direction that helps answer the question.  For instance, given two classes $\bK^0\subset \bK^1$ and modules $M \subset_{pr}N_0, N_1$ in $\bK^0$, and elements $a_\ell \in N_\ell$, we have
$$\gtp_{\bK^0}(a_0/M; N_0) = \gtp_{\bK^0}(a_1/M; N_1) \implies \gtp_{\bK^1}(a_0/M; N_0) = \gtp_{\bK^1}(a_1/M; N_1)$$
but not the reverse since the amalgams witnessing equality of Galois might not be in $\bK^0$.  This means that the maximal class $\RMod$ has the \emph{minimum} number of Galois types, but that follows straight from the definition.

Beyond Galois stability and counting types, the true goal of stability theory is the development of a nonforking/independence relation as in first-order.  There is a lot of work in Abstract Elementary Classes towards this goals (especially under the name `good $\lambda$-frames', originating in Shelah \cite{sh576}), but the work in AECs of modules is much sparser.  One standout is Mazari-Armida and Rosicky \cite{mr-rel-inj}, which uses results of Lieberman-Rosicky-Vasey \cite{lrv-fork-ind} to build stable independence in an AEC of modules $\bK$ satisfying \cite[Hypothesis 2.1]{mr-rel-inj}, which includes $(\bK, \subset_{\pr})$ that are closed under direct sums, pure submodules, and pure images.  This tool is then used to prove a purely module-theoretic result:

\begin{fact}[{\cite[Theorem 4.8]{mr-rel-inj}}]
Suppose $(\bK, \subset_{\pr})$ satisfies \cite[Hypothesis 2.1]{mr-rel-inj}.  Then $M\in \bK$ is $\bK$-pure-injective iff for every pure embedding $f:A \to B$ in $\bK$ with $\|B\| \leq |R|+\aleph_0$ and $R$-homomorphism $g:A \to M$, there is $R$-homomorphism $h:B \to M$ with $h\circ f = g$.
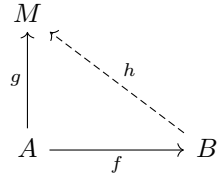
\[\begin{tikzcd}
	M \\
	\\
	A && B
	\arrow["g", from=3-1, to=1-1]
	\arrow["f"', from=3-1, to=3-3]
	\arrow["h"', dashed, from=3-3, to=1-1]
\end{tikzcd}\]

The same is true if `pure' is replaced with `RD' or `' (so just module embeddings).
\end{fact}

Another direction is undertaken by Shelah \cite{sh977}.  Shelah studies classes of modules axiomatized by the logic $\bL_{\lambda, \theta}$--allowing infinitary quantifiers and moving outside of the realm of abstract elementary classes--and proves a number of results that appear familiar to the first-order theory of modules (we suppress much of the technical detail and notation in these statements):

\begin{fact}\label{sh977-fact}
Let $\bK$ be a class of modules axiomatized in $\bL_{\lambda, \theta}$ in a `module-friendly' langauge.
\begin{enumerate}
    \item {\bf p.p. elimination of quantifiers:} For each $M \in \bK$, every formula in $\bL_{\infty, \theta}$ is equivalent to a Boolean combination of (infinitary) p. p. formulas over a `small' parameter set. (\cite[Theorem 2.4]{sh977})
    \item \label{sh977-stab-cond} {\bf Syntactic stability:} Given a collection of (infinitary) p. p. formulas $\Lambda$, then for any $M \in \bK$ of size $\lambda = \left(\lambda^{<\theta}\right)^{|\Lambda|}$, $|S_\Lambda(M)| = \lambda $. (\cite[Theorem 3.3.(1)]{sh977})
    \item {\bf Undefinability of order:} There is no formula $\phi(\bx, \by)$ in $\bL_{\infty, \theta}$ that defines a linear order of size 4 or greater on some elements from a model in $\bK$
\end{enumerate}
\end{fact}

There are additional results that suggest other syntactic aspects of stability in development.  However, these results are tempered by the fact that they are all syntactic in nature (e.g., dealing with types as sets of formulas rather than the semantic Galois types).  In first-order model theory, the compactness theorem makes this distinction immaterial.  However, since we lack compactness in Abstract Elemenatry Classes and beyond, this distinction is critical.  For instance, given two tuples $\ba \in M$ and $\bb \in N$ could have the same syntactic but different Galois types as we have no way to build the models and embeddings that would match them up.  This means that a nonelementary class of modules could have more Galois types than syntactic types, so Fact \ref{sh977-fact}.(\ref{sh977-stab-cond}) does not give an answer to Conjecture \ref{marco-quest}.

\subsection{Superstability and Uniqueness of Limit models} \label{sstab-ssec} \subsubsection{Limit models in Abstract Elementary Classes} An important object of study in Abstract Elementary Classes are \emph{limit models}.  These were introduced as `brimmed models' in \cite[Definition 1.16]{sh600}, and are especially explored in the `local study' of AECs, especially around assumptions of amalgamation and stability (see Proposition \ref{exist-limit-prop}).

\begin{defin}\label{limit-def}
    Let $\bK$ be an Abstract ELementary Class, $\lambda \geq \LS(\bK)$, and limit ordinal $\alpha <\lambda^+$.
    \begin{enumerate}
        \item $M \in \bK_\lambda$ is a \emph{$(\lambda, \alpha)$-limit model} iff there is a continuous $\prec_\bK$-increasing chain $\{ M_i \in \bK_\lambda : i< \alpha\}$ such that
        \begin{enumerate}
            \item $M_{i+1}$ is universal over $M_i$ and
            \item $M = \bigcup_{i<\alpha} M_i$
        \end{enumerate}
        \item \label{ulm-cond} $\bK$ has \emph{uniqueness of limit models in $\lambda$} iff for every $\alpha, \beta < \lambda^+$, if $M$ is a $(\lambda, \alpha)$-limit model and $N$ is a $(\lambda, \beta)$-limit model, then 
        $$M\cong N$$
    \end{enumerate}
\end{defin}

Under assumptions of amalgamation and Galois stability, one can construct universal extensions of the same size, which lets us build limit models for all parameters; (3) below shows that these assumptions are necessary.

\begin{prop}\label{exist-limit-prop}
Let $\bK$ be an Abstract Elementary Class and $\lambda \geq \LS(\bK)$.
\begin{enumerate}
    \item If $\bK_\lambda$ has amalgamation and is Galois stable in $\lambda$, then for every $M \in \bK_\lambda$, there is universal $N \in \bK_\lambda$ such that$M \prec N$.
    \item If $\bK_\lambda$ has amalgamation and is Galois stable in $\lambda$, then for every limit $\alpha < \lambda^+$, there is a $(\lambda, \alpha)$-limit model in $\bK_\lambda$.
    \item Suppose that for every $M\in \bK_\lambda$, there is universal $N \in \bK_\lambda$ such that$M \prec N$.  Then $\bK_\lambda$ has amalgamation and is Galois stable in $\lambda$.
\end{enumerate}
\end{prop}

\begin{prop}
Suppose that $M, N \in \bK$ such that $M$ is a $(\lambda, \alpha)$-limit model and $N$ is a $(\lambda, \beta)$-limit model.  If $\cf \alpha = \cf \beta$, the $M \cong N$.  If the witnessing chains have the same base model $M_0$, then the isomorphism can be built to fix $M_0$.
\end{prop}

So Uniqueness of Limit Models strengthens this proposition to all ordinals.

One might wonder why this property is a meaningful one to investigate, especially when one learns that limit models are not an object of study in elementary classes.  Indeed, \cite[Theorem 7]{gvv-ulm} seems to be the first to explore this notion in elementary classes.

The answer to this riddle is the study of when the union of saturated models is saturated.  Saturated models (see \cite[Definition 8.13]{baldwinbook}, especially $M$ is saturated iff it is $\|M\|$-saturated) are well-studied in stable theories.  In particular, for a stable theory first-order $T$, there is an important cardinal $\kappa(T)< |T|^+$ with several different meanings, including measuring the local character of the resulting independence relation (so smaller $\kappa(T)$ indicate a stronger independence relation).

Another way of understanding $\kappa(T)$ is to look at an increasing chain $\{M_i : i < \alpha\}$ of $\lambda$-saturated models and ask if the union $\bigcup_{i<\alpha} M_i$ is also $\lambda$-saturated.  A first approximation would be to compare $\lambda$ and $\alpha$: if $\lambda < \cf \alpha$, then any $<\lambda$-sized set of parameters from the union appears already in some $M_i$; thus, any type over it is realized in $M_i$ and the union is $\lambda$-saturated.  However, this first approximation could never show that a short union of highly saturated models is highly saturated.

A key result from first-order classification theory \cite[Theorem III.3.11]{sh:c} gives a better approximation that replaces $\lambda$ with $\kappa(T)$: the union of $\lambda$-saturated models is $\lambda$-saturated if the length is $\geq \kappa(T)$.  In fact this can be reversed \cite{ag-chain-sat}: if $\kappa(T) > \alpha \geq \omega$, then there is an $\alpha$-length chain of saturated models whose union is not saturated.  So the question of when unions of saturated models are saturated pins down the important cardinal $\kappa(T)$.  Note that if $\kappa(T) =\omega$, then \emph{all} unions of saturated models are saturated, and these theories are called \emph{superstable}.

Now we return to Abstract Elementary Classes and the uniqueness of limit models.  We infer from first order that the union of saturated models being saturated is an important concept to explore in nonelementary classification theory (and this has been shown true in works such as \cite{bv-unionsat}).  However, the local approach\footnote{The `local approach' and `global approach' are defined and discussed in \cite[Section 2.1]{bv-survey}.  Briefly, the local approach (exemplified by \cite{sh576}) starts with assumptions on $\bK_\lambda$ to prove properties of $\bK_{\lambda^+}$ and above, while the global approach (exemplified by \cite{sh394}) makes assumptions on all of $\bK$.} to analyzing Abstract Elementary Classes does not allow us to use saturated models! We work with Galois types over models and restricting our focus to $\bK_\lambda$ means there are no smaller models to look at types over.

The use of limit models can then be explained as the `local shadow' of Galois saturated models, and the uniqueness of limit models is the corresponding shadow of the union of saturated models being saturated.  The following proposition makes this precise.

\begin{prop}\label{local-shadow-prop}
Suppose $\bK$ is an Abstract Elementary Class and $\lambda > \LS(\bK)$ such that$\bK_\lambda$ has amalgamation and is Galois stable in $\lambda$.
\begin{enumerate}
    \item If $\alpha < \lambda^+$, then every $(\lambda, \alpha)$-limit model is $\cf \alpha$-Galois saturated if $\cf \alpha > \LS(\bK)$.
    \item If $\alpha < \lambda^+$, then every $(\lambda, \alpha)$-limit model can be written as a union of Galois saturated models from $\bK_\lambda$ of length $\alpha$.
    \item $\bK_\lambda$ satisfies the uniqueness of limit models  iff the union of Galois saturated models in $\bK_\lambda$ of length $<\lambda^+$ is Galois saturated.
\end{enumerate}
\end{prop}

More recently, \cite[Conjecture 1.1]{bv-limitstable} has begun the investigation of the set of regular $\alpha < \lambda^+$ such that $(\lambda, \alpha)$-limit models are isomorphic to $(\lambda, \lambda)$-limit models.  The conjecture is that this set is upwards closed and the minimum value in it plays a similar role for AECs that $\kappa(T)$ does for first-order.  \cite{bm-spec-lim} has continued this line.

So far, the term `superstability' has not been defined.  One reason for this is that the several equivalent characterizations of superstability in elementary classes bifurcate into distinct properties in Abstract Elementary Classes (although \cite{gv-superstable} shows that they reconverge in tame Abstract Elementary Classes).  For its most basic definition, superstability is used to mean `Galois stable on a tail of cardinals.'. However, one of the other possible definitions of superstability is the uniqueness of limit models, and some sources have begun to use this meaning as its definition.

\subsubsection{Limit models in classes of modules}

Prior to Mazari-Armida's work \cite{m-desc-limit}, the study of limit models and their uniqueness had been largely theoretical: it was explored in abstract set ups, but no examples had been explored.  \cite{m-desc-limit} examined the class of torsion-free abelian groups and showed that a $(\lambda, \alpha)$-limit model with $\cf \alpha \geq \omega_1$ must be algebraically compact \cite[Lemma 4.10]{m-desc-limit} owing to the uncountable cofinality of $\alpha$ and the countable condition on algebraic compactness, which then allows one to deduce that all such limits must have the same structure \cite[Lemma 4.14]{m-desc-limit}.  On the other hand, $(\lambda, \omega)$-limit models are not isomorphic to the other models (and thus not algebraically compact).  The argument for this relies heavily on model theoretic work: if all limit models were isomorphic, then this gives a form of superstability, namely \cite[Fact 4.19]{m-desc-limit} combines work of Grossberg and Vasey \cite{v-superstabaec,gv-superstable} to show that the nice AEC conditions that $\bK^{tf}$ satisfies plus uniqueness of limit models would imply Galois stability on a tail of cardinals; however, $\bK^{tf}$ has long been known to be only strictly stable (in particular, stable in exactly the cardinals $\lambda = \lambda^{\aleph_0}$, see \cite{bet-nperp}).\footnotei{Marcos: Surely one can construct a $(\lambda, \omega)$-limit model that is not algebraically compact by hand, right?}

Similarly, $\bK^{tor}$ is shown \cite[Lemma 5.7]{m-note-tor} to be strictly stable by showing $(\lambda, \omega)$-limit and $(\lambda, \omega_1)$-limit models are not isomorphic by giving explicit descriptions of them.

Further work on limit models picked up on a familiar theme in ring theory: determining properties of the ring $R$ by understanding properties of $R$-modules.
\begin{fact}\
\begin{enumerate}
    \item \label{noeth-cond} $R$ is left Noetherian iff $\RMod$ has uniqueness of limit models \cite[Theorem 3.12]{m-sstab-ssimp}

    \item $R$ is strictly left $<\aleph_n$-Noetherian iff $\RMod$ has exactly $n+1$ nonisomorphic limit models of size $\lambda$ for each $\lambda \geq |R|^+ + \aleph_1$ where $\Rmod$ is $\lambda$-stable. \cite[Theorem 3.17]{m-lim-noeth}.

    \item \label{pss-cond} $R$ is left pure-semisimple iff $(\RMod, \subset_\pr)$ has uniqueness of limit models. \cite[Theorem 4.28]{m-sstab-ssimp}

    \item $R$ is left perfect iff $(\Rmod^{flat}, \subset_{\pr})$ has uniqueness of limit models in some/every $\lambda \geq |R|^++\omega_1$. \cite[Theorem 3.15]{m-sstab-flat}

\end{enumerate}
\end{fact}
The proofs of these different results are different.  One neat commonality is that the proof of uniqueness of limit models often goes through a Schroeder-Bernstein style argument: for instance, in \ref{noeth-cond}) and (\ref{pss-cond}), one shows that any two limit models of $\RMod$ are injective and that they are embeddable into each other because they are universal, then uses a result of Bumby (see \cite[Fact 2.26]{m-sstab-flat}) to conclude they are isomorphic.

\bibliography{bib}
\bibliographystyle{alpha}

\end{document}